\documentclass[a4paper,11pt,leqno]{amsart}

\usepackage{amsmath}
\usepackage{amsthm}
\usepackage{amsfonts}
\usepackage{amssymb}

\usepackage{tikz}
\usetikzlibrary{3d, calc,shadings}
\usepackage{enumerate}
\usepackage{colonequals}
\usepackage{booktabs}
\usepackage{url}
\usepackage{todo}
\usepackage{subcaption}
\usepackage{anyfontsize} 
\usepackage{pdflscape}
\usepackage{algorithm}
\usepackage{algpseudocode}

\newcommand{\inv}{^{-1}}

\newcommand{\eps}{\varepsilon}
\newcommand{\N}{\mathbb N}

\newcommand{\R}{\mathbb R}
\newcommand{\rn}{\R^n}

\newcommand{\bS}{\mathbb S}

\newcommand{\vecv}{\mathbf{v}}
\newcommand{\vecw}{\mathbf{w}}
\newcommand{\vecu}{\mathbf{u}}
\newcommand{\ip}[2]{\langle #1, #2 \rangle}

\newtheorem{theorem}{Theorem}[section]
\newtheorem*{theorem*}{Theorem}{\bf}{\it}

\newtheorem*{proposition*}{Proposition}{\bf}{\it}
\newtheorem{lemma}[theorem]{Lemma}
\newtheorem*{lemma*}{Lemma}{\bf}{\it}

\theoremstyle{definition}
\newtheorem{definition}[theorem]{Definition}
\newtheorem*{definition*}{Definition}
\theoremstyle{remark}
\newtheorem{remark}[theorem]{Remark}
\numberwithin{equation}{section}

\numberwithin{equation}{section}

\tikzset{
  cube corner/.pic={
    \def\L{1}%
    \coordinate (O) at (0,0,0);
    \coordinate (X) at (\L,0,0);
    \coordinate (Y) at (0,\L,0);
    \coordinate (Z) at (0,0,\L);
    \draw[thick] (O) -- (X);
    \draw[thick] (O) -- (Y);
    \draw[thick] (O) -- (Z);
    \draw[dashed] (X) -- ($(X)+(0,0,\L)$) -- ($(Y)+(0,0,\L)$) -- (Y);
    \draw[dashed] (X) -- ($(X)+(0,\L,0)$);
    \draw[dashed] (Y) -- ($(Y)+(\L,0,0)$);
    \begin{scope}
      \clip (0,0,0) circle (0.5);
      \shade[ball color=blue!40!white] (0.15,0.15,0.1) circle (0.5);
      \draw[black] (0,0) circle (0.5);
    \end{scope}
  }
}

\title[A short survey on almost orthogonal vectors]{A short survey on almost orthogonal vectors in a few specific large dimensions}
\author{$\R$ami Luisto}
\address{
Faculty of Information Technology, P.O. Box 35 (Mattilanniemi 2), FI-40014 University of Jyv\"askyl\"a, Finland
\and
Digital Workforce Services
Mechelininkatu 1 a,
00180 Helsinki, Finland
\and
Department of Mathematics and Statistics, P.O. Box 68 (Pietari Kalmin
katu 5), FI-00014 University of Helsinki, Finland}
\email{rami.luisto@gmail.com}

\subjclass[2010]{52C17 (68T07, 60D05, 46B85, 94B65)}
\date{\today}

\begin{document}

\maketitle

\begin{abstract}
    The concept of \emph{almost orthogonal vectors}, i.e.\ vectors whose cosine similarity is close to $0$, relates to topics both in pure mathematics and in coding theory under the guises of spherical packing and spherical codes. In recent years the rise of advanced language models in AI has created new interest in this concept as the models seem to store certain concepts as almost orthogonal directions in high-dimensional spaces. In this survey we represent some ideas regarding almost orthogonal vectors through three approaches: (1) the mathematical theory of almost orthogonality, (2) some observations from the embedding spaces of language models, and (3) generation of large sets of almost orthogonal vectors by simulations.
\end{abstract}

\tableofcontents

\section{Introduction}
\label{sec:Intro}

In a given Euclidean space $\R^n$ there can exist at most $n$ directions that are pairwise orthogonal \cite{axler2015linear}. If we relax the condition of orthogonality by requiring the inner product of unit vectors to be merely \emph{close} to 0 the situation changes drastically, especially in higher dimensions. This phase transition in high and low dimensions relates to a concept sometimes known as \emph{the curse (or blessing) of high dimensionality}. We'll discuss this further in the coming sections, see especially Section \ref{sec:highandlow} 

We say that two vectors $\vecv, \vecw \in \rn$ are \emph{$\eps$-almost orthogonal} for $\eps > 0$ if their cosine similarity lies between $-\eps$ and $\eps$. The main motivating question for this paper is as follows:
\begin{quote}
    How many pairwise $\eps$-almost orthogonal vectors can fit in $\rn$ for various values of $\eps$ and $n$? 
\end{quote}
As we'll discuss in Section \ref{sec:sphericalcodesandspherepacking} this question is closely tied to various open problems related to both \emph{spherical codes} and \emph{sphere packing}. In particular, exact solutions to the best configurations of vectors are beyond the scope of this survey, and instead we will focus on various estimates.

Since e.g.\ spherical codes have many practical applications, there already exists literature on approximate solutions. However, these do tend to focus more on low-dimensional examples -- see e.g.\ \cite{cohn2023spherical} where the examples only go up to dimension 32. This is natural if you aim to use this with e.g.\ binary blocks that tend to have word sizes of 8, 16, or 32 bits. By contrast, the existing results for almost orthogonal vectors in higher dimensions tend to focus on the \emph{asymptotic} behavior of the results, i.e.\ what happens when the ambient dimension grows very large. In particular, the typical results like the Johnson-Lindenstrauss Lemma tend not to yield good estimates in dimensions below several thousands -- see Section \ref{sec:JohnsonLindenstrauss}.

The contextual motivation for us, on the other hand, arises from the so called \emph{embedding vectors} of various AI-models. For a contemporary example, various language models based on the transformer architecture (including GPT, LLaMa, Gemini and Mistral) encode the meaning of text through embedding vectors that live in high-dimensional Euclidean spaces like $\R^{768}$ or $\R^{1024}$ (see e.g.\cite{ethayarajh2019contextual, mistral2023mistral, rogers2021primer}). These dimensions lie between the small dimensions of coding theory and high dimensions of asymptotical results. This gap is of primary interest in this survey.

Based on current understanding, see e.g.\ \cite{elhage2022toy,henighan2023superposition,bricken2023towards,templeton2024scaling}, modern language models models seem to store different semantic concepts as different directions in these embedding spaces. Crucially for us, these different concepts seem to be encoded as \emph{almost} orthogonal vectors, and it turns out that the geometry of e.g.\ $\R^{768}$ can encode much more such directions than just 768, though lossily. A driving motivation in this survey is to better understand how many almost orthogonal vectors can be fitted in various high-dimensional Euclidean spaces where the particular dimensions we study are the embedding space dimensions of particular contemporary language models -- for the sake of exposition we will often focus on dimensions 512, 768, 1024 and 2048 when concrete examples are needed.

As per the parameter $\eps$ in $\eps$-almost orthogonality we will in most examples focus on either $\eps=0.1$ or $\eps=0.01$ when we need to fix the parameter ourselves. This is again purely for the sake of exposition, as the values themselves do not carry any special meaning. Note that $\eps=0.1$ corresponds to angles between 84 and 96 degrees, while $\eps=0.01$ corresponds to angles of roughly $90\pm 0.5$ degrees. Also note that if unit vectors $\vecv$ and $\vecw$ have an inner product of $0.1$, then $\vecv = 0.1 \vecw + 0.9\vecu$ where $\vecu$ is orthogonal to $\vecw$, meaning that "10\% of $\vecv$ is explained by $\vecw$".

Part of the codebase used for the simulations in this work is available at \url{https://github.com/ramiluisto/AlmostOrthogonalGenerators}.

\subsection{Structure of this survey}

We'll start by "setting the scene" and describe in a bit more detail on why we are interested in almost orthogonality in high-dimensional spaces when studying modern (Large) Language models. We will then go through some mathematical preliminaries, including some more heuristical comments on how the geometry of low and high dimensions differs. We'll then move on to studying some (trivial) mathematical bounds we can get to the amount of almost orthogonal vectors in Euclidean spaces. We then look at some other classical results, in particular the Johnson-Lindenstrauss Lemma. The bounds we see here turn out to be somewhat meagre compared to practical situations. Thus we turn into simulations and compare various methods of generating almost orthogonal vectors. In the final section we'll draw conclusions on all of this. 

We wish to emphasize that this paper is structured as a \emph{survey}; it claims no new mathematical theorems, rather it tries to provide a convenient summary of useful known information. On the other hand it does not claim to be a \emph{comprehensive} survey as many areas are only briefly mentioned, and we lack the expertise to comment authoratively to many topics discussed here.

\textbf{Acknowledgements:} This work started as an appendix to a master's thesis. We would like to thank the advisors Sami \"Ayram\"o and Ida Toivanen for their support. We furthermore extend our thanks to Ilkka P\"ol\"onen for discussions on the topic.

\section{Context -- Latent spaces of (Large) Language Models}

We will not repeat here a detailed introduction to the architecture of transformer models, rather we given a quick overview of the terminology. For a proper introduction, see e.g.\ \cite{tunstall2022natural}.

Most classical transformer-based language models break incoming text into \emph{tokens}. These tokens tend to be common words or word fragments. The collection of tokens a model has is called the \emph{vocabulary}. The vocabulary size of e.g.\ the different BERT\footnote{BERT is what I would call the archetypical small language model, see e.g.\ \cite{rogers2021primer} for a comprehensive text on the topic.} variants are in the ballpark of 30--50 thousand tokens. For a more thorough introduction to tokens and their usage we refer to \cite{devlin2018bert} and \cite{mikolov2013efficient}.

For each token in the vocabulary the model has a learned embedding vector of a particular dimension $d_\text{model}$. These vectors are called the \emph{input embeddings} of the model. Thus a given input sequence of natural text is broken into a sequence of tokens, and each of these tokens is then converted into the corresponding embedding vector. This sequence of embedding vectors is then passed through several \emph{transformer blocks} by the model. These transformer blocks update each of the embedding vectors of the input sequence based on other vectors in the sequence. This evolution of the sequence of embedding vectors is sometimes called \emph{the residual stream}. For the standard transformer architecture, the $d_\text{model}$ stays the same for the embeddings as they progress through the transformer blocks.\footnote{Though e.g.\ the ALBERT model uses factorized embedding; token embeddings are 128-dimensional but are expanded internally to a 768-dimensional vector at various points in the processing.} We've listed some of the values\footnote{We list the exact dimensions when known, but for some closed source models list the best known guesses based on public information available.} of $d_\text{model}$ in Tables \ref{tab:embedding_space_dimensions_encoders}, \ref{tab:embedding_space_dimensions_decoders} and \ref{tab:embedding_space_dimensions_encoder_decoder} for a few common SLM and LLM architectures.

\begin{table}[h!]
    \centering
    \small
    \begin{tabular}{lr}
    \toprule
    \textbf{Model} & $d_\text{model}$ \\
    \midrule
    \multicolumn{2}{l}{\textit{Encoder-only (BERT-style)}} \\
    BERT (base) \cite{devlin2018bert} & 768 \\
    BERT (large) \cite{devlin2018bert} & 1024 \\
    RoBERTa (base) \cite{liu2019roberta} & 768 \\
    RoBERTa (large) \cite{liu2019roberta} & 1024 \\
    DistilBERT \cite{sanh2019distilbert} & 768 \\
    ALBERT (base) \cite{lan2019albert} & 128 \\
    XLM-RoBERTa \cite{conneau2019xlmr} & 768 \\
    XLNet (base) \cite{yang2019xlnet} & 768 \\
    XLNet (large) \cite{yang2019xlnet} & 1024 \\
    ERNIE \cite{sun2019ernie} & 768 \\
    ModernBERT (base) \cite{warner2024smarter} & 768 \\
    ModernBERT (large) \cite{warner2024smarter} & 1024 \\
    \bottomrule
    \end{tabular}
    \caption{Token-embedding (input) vector sizes for a wide range of language model families of the \emph{encoder} architecture.}
    \label{tab:embedding_space_dimensions_encoders}
\end{table}

\begin{table}[h!]
    \centering
    \small
    \begin{tabular}{lr}
    \toprule
    \textbf{Model} & $d_\text{model}$ \\
    \midrule
    \multicolumn{2}{l}{\textit{Decoder-only (GPT, LLaMA, etc.)}} \\
    GPT-1 \cite{radford2018improving} & 768 \\
    GPT-2 (small) \cite{radford2019gpt2} & 768 \\
    GPT-2 (medium) \cite{radford2019gpt2} & 1024 \\
    GPT-2 (large) \cite{radford2019gpt2} & 1280 \\
    GPT-2 (XL) \cite{radford2019gpt2} & 1600 \\
    GPT-3 (175B) \cite{brown2020gpt3} & 12288 \\
    GPT-3.5 & $\sim$12288$^\dagger$ \\
    GPT-4 & Not disclosed \\
    Codex (12B) & $\sim$5120$^\dagger$ \\
    LLaMA-1 (7B) \cite{touvron2023llama} & 4096 \\
    LLaMA-1 (13B) \cite{touvron2023llama} & 5120 \\
    LLaMA-1 (33B) \cite{touvron2023llama} & 6656 \\
    LLaMA-1 (65B) \cite{touvron2023llama} & 8192 \\
    LLaMA-2 (7B) \cite{touvron2023llama2} & 4096 \\
    LLaMA-2 (13B) \cite{touvron2023llama2} & 5120 \\
    LLaMA-2 (70B) \cite{touvron2023llama2} & 8192 \\
    Mistral (7B) \cite{mistral2023mistral} & 4096 \\
    Mixtral (8$\times$7B) \cite{jiang2024mixtral} & 4096$^\dagger$ \\
    Bloom (176B) \cite{le2023bloom} & 14336 \\
    Falcon (40B) \cite{almazrouei2023falcon} & $\sim$8192$^\dagger$ \\
    Grok-1 \cite{xai2024grok} & Not disclosed \\
    \bottomrule
    \end{tabular}
    \caption{Token-embedding (input) vector sizes for a wide range of language-model families of the \emph{decoder} architecture. $^\dagger$Estimated values; "Not disclosed" indicates the vendor has not released architecture details.}
    \label{tab:embedding_space_dimensions_decoders}
\end{table}

\begin{table}[h!]
    \centering
    \small
    \begin{tabular}{lr}
    \toprule
    \textbf{Model} & $d_\text{model}$ \\
    \midrule
    \multicolumn{2}{l}{\textit{Encoder–decoder (T5, PaLM, …)}} \\
    T5 (small) \cite{raffel2019t5} & 512 \\
    T5 (base) \cite{raffel2019t5} & 768 \\
    T5 (large) \cite{raffel2019t5} & 1024 \\
    T5 (3B) \cite{raffel2019t5} & 1024 \\
    T5 (11B) \cite{raffel2019t5} & 1024 \\
    PaLM (540B) \cite{chowdhery2023palm} & 18432 \\
    Gopher (280B) \cite{rae2021scaling} & 16384 \\
    Chinchilla (70B) \cite{hoffmann2022training} & 8192 \\
    AlphaCode (41B) \cite{li2022competition} & 6144 \\
    Claude 1/2 \cite{anthropic2023claude2} & Not disclosed \\
    Gemini Ultra \cite{google2023gemini} & Not disclosed \\
    \bottomrule
    \end{tabular}
    \caption{Token-embedding (input) vector sizes for a wide range of language-model families for the \emph{encoder-decoder} architecture. "Not disclosed" indicates the vendor has not released architecture details.}
    \label{tab:embedding_space_dimensions_encoder_decoder}
\end{table}

For us the input embeddings are of particular interest because they are context free -- the embeddings do not depend on any way on the other tokens in the sequence nor on the relative or absolute position of a particular token. Thus we can study these vectors as "independent embedding vectors". For latter parts in the residual stream any embedding vector is affected by the other tokens and their positions in the input -- i.e.\ the context. Thus the study of the embedding vectors and their geometry in the latter layers of transformers usually requires sampling of text examples and their careful analysis. The selection of texts to use naturally has a considerable effect, so for the examples in this section we simply restrict ourselves to the input embeddings of various models.

\begin{figure}[ht]
    \centering
    \includegraphics[width=1.0\textwidth]{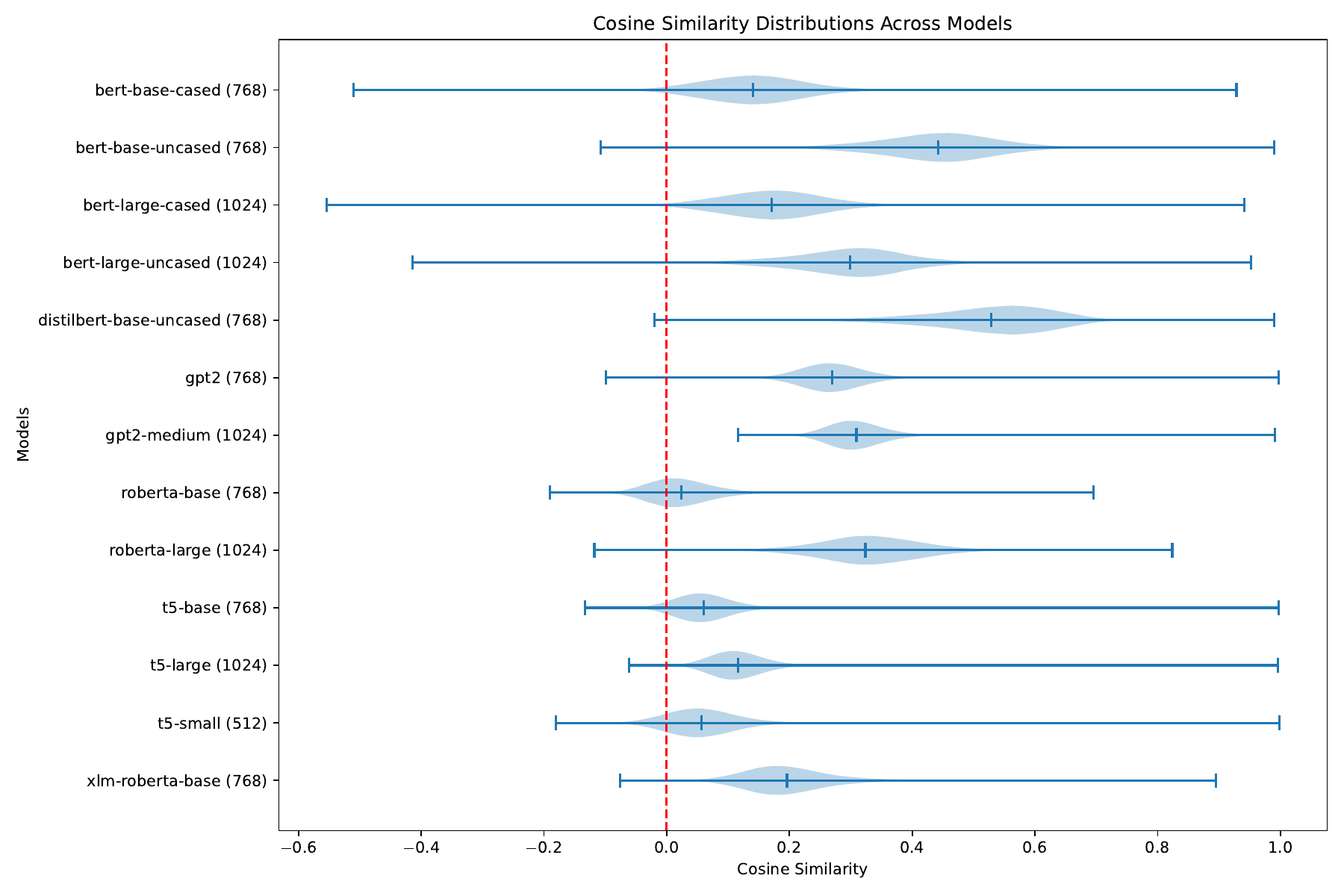}
    \caption{The distribution of pairwise cosine similarities of input embeddings for various models. The dimension of the embedding space is in parenthesis after the model name.}
    \label{fig:violin_plot_for_model_cosims}
\end{figure}

To get an idea on what kind of cosine similarities we see in the real world, we've taken a few common transformer-based models, extracted their input embeddings and calculated the cosine similarities between each embedding vector pair (excluding the trivial cosine similarities of a vector with itself). In Figure \ref{fig:violin_plot_for_model_cosims} we plot the distribution of these cosine similarities for each model.

A noteworthy property noticed already in \cite{mu2017all} and \cite{ethayarajh2019contextual} is that the cosine similarities are not at all equally distributed around the origin. Instead for almost all models the main mass of cosine similarities is strictly positive, sometimes by a wide margin. With the interpretation that different directions encode various semantic meanings in any transformer model, the conclusion here would then be that almost all transformer models tend to see quite a lot of similarities between any two words or tokens -- at least at the input embedding phase of things. Even if there are some meanings that are somewhat orthogonal, strong negative correlation in the form of a cosine similarity near $-1$ seems to be unheard of.
We note that it is not the "purpose" of the language models or their input embedding vectors to minimize pairwise absolute cosine similarities. The models do want to understand differences between tokens, but they also really "want" to find similarities. Furthermore even if the model uses some "basis" collection of almost orthogonal directions to encode various meanings that correlate more or less, most of the embedding vectors are then in a sense mapped to be linear combinations of these major dimensions and not dominant directions themselves.

We've listed some more detailed statistical numbers in Tables \ref{table:model_stats_cosine} and \ref{table:model_stats_norms}.

\begin{table}[!ht]
\centering
\tiny
\begin{tabular}{lrrrrrr}
\toprule
Model & Emb.\ Dim & Mean CosSim & Std CosSim & CosSim 25\% & CosSim 50\% & CosSim 75\% \\
\midrule
xlm-roberta-base & 768 & 0.1959 & 0.0695 & 0.1487 & 0.1896 & 0.2358 \\
t5-small & 512 & 0.0569 & 0.0592 & 0.0170 & 0.0537 & 0.0925 \\
t5-large & 1024 & 0.1161 & 0.0472 & 0.0849 & 0.1126 & 0.1425 \\
t5-base & 768 & 0.0607 & 0.0491 & 0.0281 & 0.0576 & 0.0889 \\
roberta-large & 1024 & 0.3237 & 0.0912 & 0.2727 & 0.3273 & 0.3823 \\
roberta-base & 768 & 0.0234 & 0.0543 & -0.0136 & 0.0189 & 0.0550 \\
gpt2-medium & 1024 & 0.3091 & 0.0471 & 0.2775 & 0.3058 & 0.3365 \\
gpt2 & 768 & 0.2697 & 0.0537 & 0.2349 & 0.2679 & 0.3022 \\
distilbert-base-uncased & 768 & 0.5283 & 0.1023 & 0.4668 & 0.5398 & 0.5998 \\
bert-large-uncased & 1024 & 0.2986 & 0.0999 & 0.2397 & 0.3036 & 0.3598 \\
bert-large-cased & 1024 & 0.1713 & 0.0859 & 0.1155 & 0.1716 & 0.2261 \\
bert-base-uncased & 768 & 0.4424 & 0.0983 & 0.3811 & 0.4453 & 0.5052 \\
bert-base-cased & 768 & 0.1406 & 0.0848 & 0.0845 & 0.1397 & 0.1944 \\
\bottomrule
\end{tabular}
\caption{Cosine similarity statistics: Embedding dimension, mean/std and quartiles of the cosine similarity.}
\label{table:model_stats_cosine}
\end{table}

\begin{table}[!ht]
\centering
\tiny
\begin{tabular}{lrrrrrr}
\toprule
Model & Emb.\ Dim & Mean Norm & Std Norm & Norm 25\% & Norm 50\% & Norm 75\% \\
\midrule
xlm-roberta-base & 768 & 5.8701 & 0.4551 & 5.6919 & 5.9038 & 6.1314 \\
t5-small & 512 & 522.4708 & 65.9444 & 483.9047 & 520.5771 & 556.7801 \\
t5-large & 1024 & 457.1157 & 63.6800 & 416.2583 & 449.0090 & 488.0606 \\
t5-base & 768 & 517.7155 & 72.0042 & 473.5927 & 510.2301 & 555.2615 \\
roberta-large & 1024 & 4.3211 & 0.5238 & 4.0697 & 4.4603 & 4.6806 \\
roberta-base & 768 & 3.6286 & 0.3494 & 3.4261 & 3.6920 & 3.8746 \\
gpt2-medium & 1024 & 3.6835 & 0.4126 & 3.4084 & 3.6646 & 3.9646 \\
gpt2 & 768 & 3.9556 & 0.4379 & 3.6681 & 3.9405 & 4.2474 \\
distilbert-base-uncased & 768 & 1.6640 & 0.2713 & 1.4537 & 1.6626 & 1.8302 \\
bert-large-uncased & 1024 & 1.4575 & 0.2002 & 1.3131 & 1.4350 & 1.5864 \\
bert-large-cased & 1024 & 1.5238 & 0.1908 & 1.3917 & 1.5058 & 1.6478 \\
bert-base-uncased & 768 & 1.4036 & 0.1989 & 1.2556 & 1.4006 & 1.5312 \\
bert-base-cased & 768 & 1.2871 & 0.1511 & 1.1886 & 1.2687 & 1.3790 \\
\bottomrule
\end{tabular}
\caption{Embedding norm statistics: Embedding dimension, mean/std and quartiles of the norm.}
\label{table:model_stats_norms}
\end{table}

\begin{figure}[ht]
    \centering
    \includegraphics[width=1.0\textwidth]{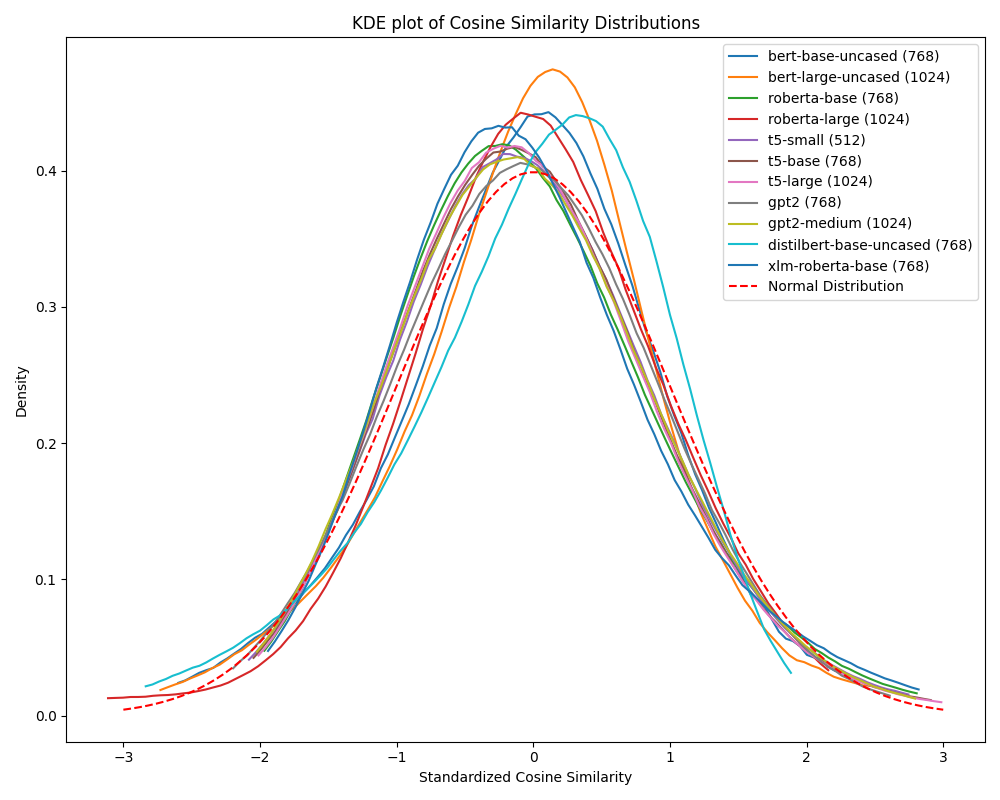}
    \caption{The distribution of pairwise cosine similarities of input embeddings for various models after being scaled as a probability distribution and normalized. Normal distribution is shown for comparison.}
    \label{fig:kde_plot_for_model_cosims}
\end{figure}

We further note that if we rescale the distributions as probability distributions and then normalize them by translating by the mean and scaling by the standard deviation, the distributions are very alike to each other and the normal distribution. See Figure \ref{fig:kde_plot_for_model_cosims}. We are not quite sure if this is expected or surprising. 

\section{Mathematical preliminaries}

We'll start off with some mathematical preliminaries. We go through basic notations and then define some special functions in order to approximate the volumes and ares of high-dimensional balls, spheres and their subsets. We then present some high-level observations on the geometrical differences of high and low dimensions.

\subsection{Notation}

We denote by $\N = \{0, 1, 2, \ldots \}$ the set of natural numbers and use the shorthand of $\N_+$ for the strictly positive natural numbers $\{ 1, 2, 3, \ldots\}$. By $\R$ we denote the set of real numbers, while $\R_+$ stands for the non-negative real numbers $[0, \infty)$.

By $\rn$ we denote the $n$-dimensional Euclidean space for $n \geq 1$, i.e. the set
\begin{align*}
    \left\{ (x_1, x_2, \ldots, x_n) \mid x_i \in \R \right\}.
\end{align*}
We often denote the vectors in $\rn$ by $\vecv$, $\vecw$, $\vecu$ or some other bold font letter. In these cases we implicitly assume that the coordinates of the vector are denoted by subscripts of an unboldened letter, e.g.\ the components of the vector $\vecw$ are $(w_1, w_2, \ldots, w_n)$, and we tacitly assume that the dimension $n$ is clear from the context. The \emph{origin} is the point with all zero coordinates and we denote it as $\mathbf{0} \colonequals (0, \ldots, 0)$.

We equip $\rn$ with the inner product
\begin{align*}
    \ip{\cdot}{\cdot} \colon \rn \times \rn \to \R, \quad
    \ip{\vecv}{\vecw} = \sum_{i=1}^n v_i w_i
\end{align*}
and the Euclidean norm
\begin{align*}
    \| \cdot \| \colon \R^n \to \R_+, \quad
    \|\vecv\| = \sqrt{\ip{\vecv}{\vecv}} = \sqrt{\sum_{i=1}^n v_i^2}
    \quad \text{ for } \vecv \in \rn.
\end{align*}
A vector with length 1 is called a \emph{unit vector}.
We call two non-zero vectors $\vecv$ and $\vecw$ \emph{orthogonal} if $\ip{\vecv}{\vecw} = 0$, and \emph{orthonormal} if they are orthogonal unit vectors.



With the inner product we can also define a crucial concept in this survey - \emph{the cosine similarity}.
\begin{definition}
    Let $\vecv, \vecw$ be two vectors in $\rn \setminus \{ \mathbf{0}\}$. We set their \emph{cosine similarity}
    to be the number
    \begin{align*}
        \operatorname{cs}(\vecv, \vecw)
        = \frac{\ip{\vecv}{\vecw}}{\|\vecv\| \cdot \|\vecw\|}.
    \end{align*}
\end{definition}
Note that orthogonal vectors have cosine similarity of $0$, and the cosine similarity of a vector with itself is always $1$. Furthermore, scaling of either vector by a non-zero scalar has no effect on the vectors' cosine similarity. 

For any two non-zero vectors $\vecv, \vecw$,
\begin{align*}
    \cos(\theta) 
    = \frac{\ip{\vecv}{\vecw}}{\|\vecv\| \cdot \|\vecw\|}
    = \operatorname{cs}(\vecv, \vecw),
\end{align*}
where $\theta$ is the angle between the vectors.

\subsubsection{A few special functions}

When studying the volumes and areas of high-dimensional balls and spheres the formulas tend to involve a few \emph{special functions}. These are akin to the so called \emph{special constants} like $\pi$ or $e$, except that the Gamma and Beta functions are \emph{functions} and not constants. They are natural objects that pop up\footnote{One could take the Platonic view that these are actually existing objects that have some special property, or a more Kolmogorov-complexity style of approach that these happen to be concepts that we come across so frequently that we've given them special short names.} in many different settings, just like $\pi$.  Regardless, they are needed for many volume and area estimates of high-dimensional balls and spheres; for further details and proofs we refer the reader to \cite{abramowitz1948handbook}.

The Gamma function can be defined in a subdomain of the complex plane, but for us it is sufficient to study it for positive real numbers.
\begin{definition}
    The \emph{Gamma function} $\Gamma$ is defined as
    \begin{align*}
        \Gamma \colon (0, \infty) \to \R, \quad
        \Gamma(x) = \int_0^\infty t^{x-1} e^{-x}\text{ d}t.
    \end{align*}
\end{definition}

For our purposes, the following properties of the Gamma function will be sufficient, see again \cite{abramowitz1948handbook} for the proofs.
\begin{enumerate}[($\Gamma$1)]
    \item The Gamma function generalizes the factorial, i.e.\ for all $n \in \N_+$
    $$
    \Gamma(n) = (n-1)! = (n-1)\cdot(n-2)\cdot\ldots \cdot 3 \cdot 2 \cdot 1.
    $$

    \item More generally, we have the following recursive relation for all $x\in (0,\infty)$
    $$\Gamma(x+1) = x\cdot \Gamma(x).$$

    \item We have an explicit value for $\frac{1}{2}$, namely 
    $\Gamma(\frac{1}{2}) = \sqrt{\pi}$.

    \item Stirling's formula approximation\footnote{We'll ignore here the exact accuracy of the approximation here. In the range of values we study it is less than a percent and, as we shall see, our estimations are not dependent on so fine values.} applies for the Gamma function, i.e.
    $$
        \Gamma(x+1) \approx \sqrt{2\pi x}\left(\frac{x}{e}\right)^x.
    $$
\end{enumerate}

\begin{definition}
    We define the \emph{incomplete Beta function} $B \colon [0,1]\times(0,\infty)^2 \to \R$ by
    \begin{align*}
        B(x;\,a,b) 
        = \int_0^x t^{a-1}\,(1-t)^{b-1}\,dt.
    \end{align*}

    And the \emph{Beta function} is defined simply by setting the parameter $x = 1$ and denoted with a slight abuse of notation as $B(a,b)$.
\end{definition}

\begin{definition}
    The \emph{regularized incomplete Beta function}
    $I \colon [0,1]\times(0,\infty)^2 \to \R$ is defined as
    \begin{align*}
        I_x(a,b)
        = \frac{B(x;\,a,b)}{B(a,b)}.
    \end{align*}
\end{definition}

\subsection{Volumes and areas of spheres and balls}\label{sec:volumes_and_areas_of_spheres_and_balls}

With the Gamma function and the various Beta function variants, we can now move on to studying the volumes and areas for spherical objects in high dimensions.

By \emph{the unit cube} in $\rn$ we mean the set $[0,1]^n$, and by \emph{the unit cube centered at the origin} the set $\left[-\frac{1}{2}, \frac{1}{2}\right]^n$.
When we talk of a \emph{cube of side length $a$} we mean a unit cube that has been scaled by the scalar $a$ in all coordinates.

In $\rn$, the volume of a cube of side length $a$ is simply $a^n$, and the diameter of such a cube is $a \sqrt{n}$. In particular the unit cube has volume $1$ and diameter $\sqrt{n}$.

We denote by $\overline{B}_n(r)$ the \emph{closed ball of radius $r$ in $\rn$}, i.e.\ the set
\begin{align*}
    \overline{B}_n(r) 
    = \{ \vecv \in \rn \mid \| \vecv \| \leq r\}.
\end{align*}
Furthermore we denote its volume by $V_n(r)$.

Similarly the $n-1$ dimensional \emph{sphere of radius $r$} is the set
\begin{align*}
    \bS
    = \{ \vecv \in \rn \mid \| \vecv \| = r\}.
\end{align*}

The volume of a ball of radius $r$ in $\R^n$ is given by
\begin{align}
    V_n(r) \colonequals \frac{\pi^{n/2}}{\Gamma(n/2 + 1)}f^n.
\end{align}
(See again \cite{abramowitz1948handbook} for details.)

With the Stirling formula we can get an approximation for the volume as
\begin{align*}
    V_n (r)
    &= \frac{\pi^{n/2}}{\Gamma(n/2 + 1)} r^n
    \approx \frac{\pi^{n/2}}{\sqrt{2 \pi (n/2 + 1)} \cdot \left(\frac{n/2 + 1}{e}\right)^{n/2+1}} r^n \\
    &= \frac{1}{\sqrt{n\pi}}\left(\frac{2\pi e}{n}\right)^{n/2}r^n.
\end{align*}
From this we can then calculate an approximate value for the function $R_n \colon \R_+ \to \R_+$ defined by setting $R_n(V)$ to be the radius $r_V$ for which $V_n(r_V) = V$:
\begin{align*}
    R_n(V) 
    \approx (\pi n)^{1/(2n)}\sqrt{\frac{n}{2\pi e}} V^{1/n}.
\end{align*}

Note also that for the unit ball we have 
\begin{align*}
    V_n(1) 
    = \frac{1}{\sqrt{n\pi}}\left(\frac{2\pi e}{n}\right)^{n/2}.
\end{align*}

\begin{figure}[ht]\centering
    \centering
    \includegraphics[width=1.0\textwidth]{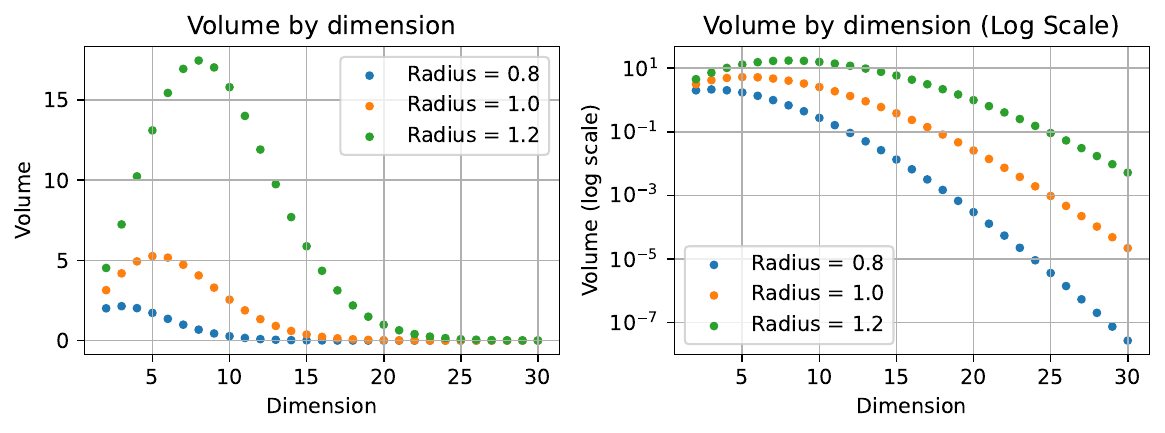}
    \caption{Volume comparisons of balls in different dimensions. The maximum value for the unit ball occurs at $n=5$.}
    \label{fig:volumeComparisons}
\end{figure}

From these formulae the crucial thing to note is that in $V_n(r)$ for large enough $n$ the value is completely dominated by the term $n^{-n/2}$, which goes to zero faster than any polynomial or exponential function. To demonstrate this, in Figure  \ref{fig:volumeComparisons} we plot some volumes of balls in various dimensions. The key takeaway here is that in higher dimensions the volume of the unit ball starts to converge to zero as $\mathcal{O}(n^{-n/2})$, i.e.\ exceedingly fast.

In dimension 768, which is important for us, trying to directly calculate the volume with the Gamma function creates an integer overflow in Python3, and probably also in most other non-symbolic math libraries. With the Stirling formula we can approximate that the unit ball in dimension 768 has volume below $10^{-300}$. Conversely, to get a ball with volume 1 we can calculate that we would need a radius of around $R_{768}(1) \approx 6.7$. This is our first important geometric observation concerning high-dimensional geometry -- \emph{the volumes of $n$-dimensional unit balls are very very small}. 

\begin{figure}[ht]\centering
\includegraphics[width=0.9\textwidth]{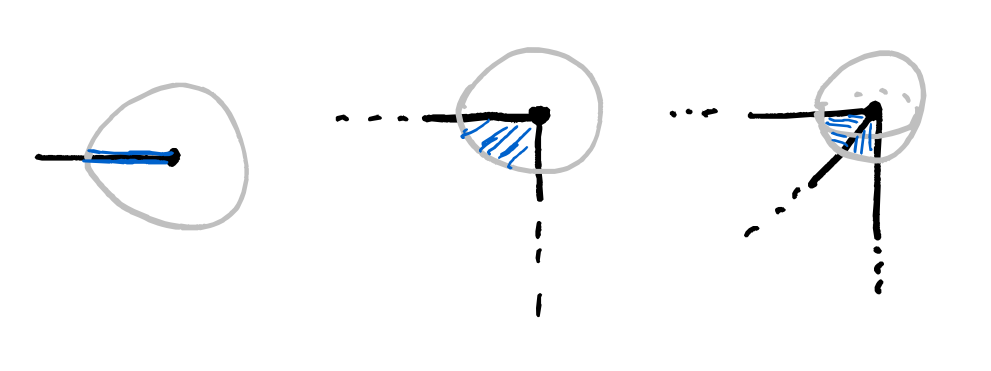}
\caption{An illustration on how much of the ambient space of a cube edge is part of the cube.}
\label{fig:cube_edge_ambient_illustration}
\end{figure}

On the other hand, as we mentioned earlier, the volume of a unit cube in $\R^n$ is $1$, and the diameter is $\sqrt{n}$. So for example in dimension 768 we get the diameter of the unit cube to be about $27.7$. This property is a first hint at a phenomenon that is sometimes described as "pointyness of high-dimensional cubes". The idea being that even though in all dimensions some parts of a unit cube centered at the origin are inside the circle (e.g.\ any point in the cube of the form $(0, \ldots, 0, 1/2, 0, \ldots, 0)$ is on the boundary of the cube and within a unit ball), the corners of the cube start to point out of the ball more and more.
Also note that in 1D, looking at the edge vertex of a cube we see that
half of the ambient space is in the cube. In 2D, only a fourth and in 3D only an eighth -- see Figure \ref{fig:cube_edge_ambient_illustration}. This theme continues in higher dimensions, with the edge vertex of a cube in dimension 768 having less than one part in $10^{-200}$ of the ambient space being part of the cube. This can be thought as another measure of the pointyness.
In dimension 768, considering that the unit cube has a diameter of around 27, we note that the corners also point out of the ball quite a lot. Furthermore there are $2^{768}$ of them, and the distance between two neighboring ones is only 1, so the cube does seem to have a more "hedgehog-style" quality to it.

\begin{figure}[ht]\centering
    \includegraphics[width=0.9\textwidth]{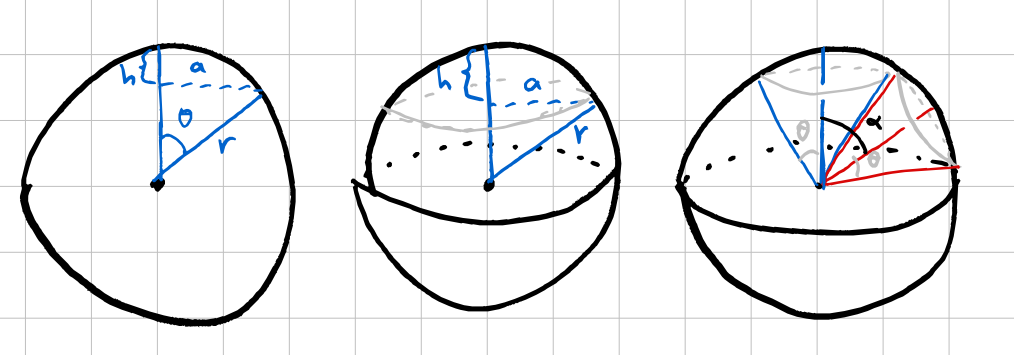}
    \caption{Example image of spherical caps.}
    \label{fig:Spherical_cap_diagram}
\end{figure}

For our geometrical considerations, the concept of \emph{spherical caps} is crucial. To construct a spherical cap on a sphere, we simply fix a point on the sphere take the intersection of the sphere with a ball centered on that point. There are various ways to parametrize the spherical cap, we refer to Figure \ref{fig:Spherical_cap_diagram} for the notation. 
In studying almost orthonormal vectors, we note that two unit vectors $\vecv$ and $\vecw$ are $\eps$-almost orthonormal exactly when on the unit sphere the spherical caps with $\vartheta = \eps/2$ around them are disjoint. 

We are interested in calculating the area\footnote{For simplicity, we'll use the term "area" even though we are not working with 2-dimensional objects. Here "$n-1$ volume" or "$n-1$ dimensional Hausdorff measure" would be more appropriate, but a mouthful.}  of such spherical caps on unit spheres, which we denote by $SC_\vartheta$, in order to calculate an area-based bound on the amount of disjoint spherical caps on a unit sphere. 

In dimension $n$ we get the following formula for the area of a spherical cap, see \cite{ledoux2001concentration} for details:
\begin{align}\label{eq:spherical_cap_area}
    A_h
    =\frac{1}{2}\frac{2\pi^{n/2}}{\Gamma\left(\frac{n}{2}\right)} r^{n-1} I_{(2rh-h^2)/r^2} \left(\frac{n-1}{2}, \frac{1}{2} \right), \quad
    0\le h\le r,
\end{align}
where $I_x$ is the regularized incomplete Beta function previously defined. For a unit cube this simplifies to 
\begin{align}
    A_h
    =\frac{1}{2}\frac{2\pi^{n/2}}{\Gamma\left(\frac{n}{2}\right)} I_{(2h-h^2)} \left(\frac{n-1}{2}, \frac{1}{2} \right), \quad
    0\le h\le 1.
\end{align}

\begin{figure}[ht]\centering
    \includegraphics[width=0.9\textwidth]{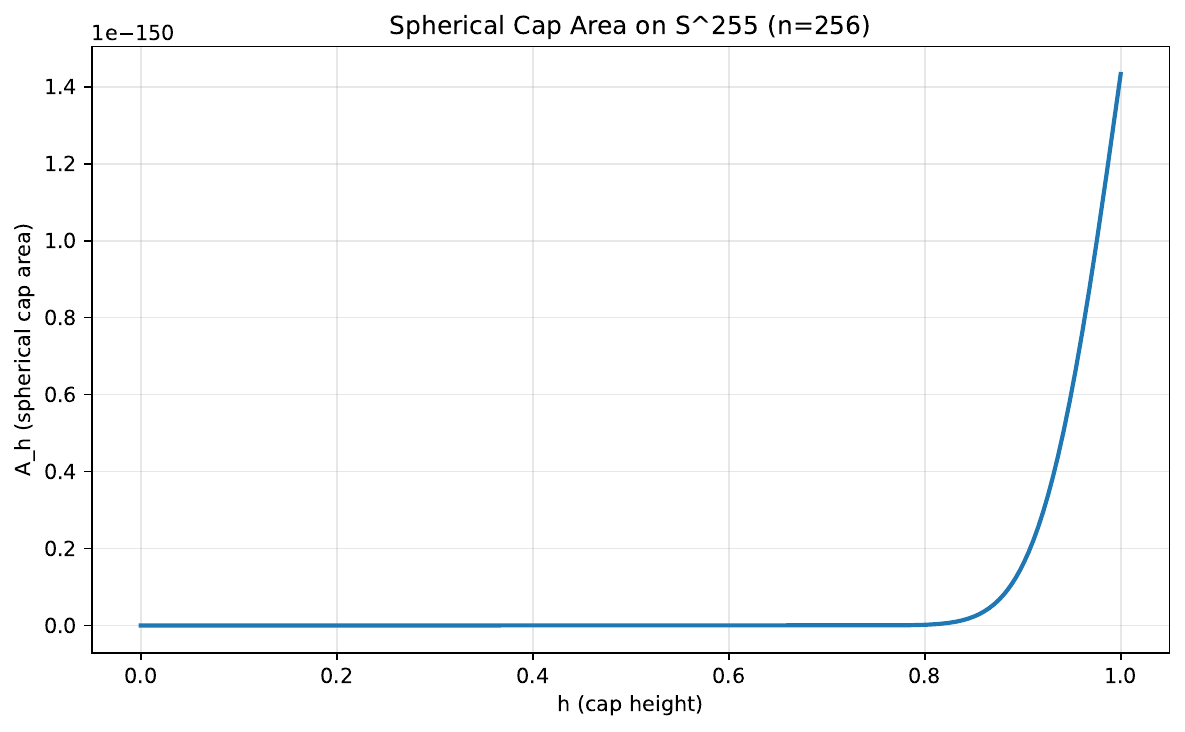}
    \caption{A plot of the area of a spherical cap of a unit sphere in $\R^{256}$.}
    \label{fig:sphere_dense_equators}
\end{figure}

We note here that if we plot $A_h$ for $h \in [0,1]$ in $\R^{256}$, see Figure \ref{fig:sphere_dense_equators}, we see that the supermajority of the area of the sphere is concentrated on the equator.

\subsection{High and low dimensions in general}
\label{sec:highandlow}

Low dimensions have several "exotic" properties, both geometrical and topological. For example:
\begin{enumerate}
    \item We've already noted above some aspects of high-dimensional geometry that are different from the more "usual" low-dimensional one. Another point of view is that it is the low dimensions that are the strange ones. In high dimensions cubes are pointy, unit balls have low volumes and even what they have is concentrated near the boundary -- indeed with the volume of a ball of radius $r$ behaving as $\mathcal{O}(r^n)$, it is easy to see that for larger $n$ the volume of $B(\mathbf{0}, 0.99)$ is considerably smaller than the volume of $B(\mathbf{0}, 1)$, meaning that the volume of the ball is concentrated near the boundary. Furthermore as we observed in Figure \ref{fig:sphere_dense_equators}, for a sphere in high dimensions the supermajority of the volume of the sphere is concentrated near the equator.\footnote{This effect is also known as "the concentration of measure", see \cite{ledoux2001concentration,vershynin2018high}.}

    \item The only way to equip a Euclidean space with a product that produces a division algebra is restricted to dimensions $1$, $2$, $4$ and $8$ \cite{adams1960non}. As a more simple example, the only way to equip a Euclidean space with a product that produces an algebraic field structure is to take $\R$ with the standard product or $\R^2$ with the product of complex numbers, see \cite{hatcher2005algebraic}.
    
    \item Knot theory works only in dimension 3; in lower dimensions a string cannot be knotted and in higher dimensions any knot can be trivially untied \cite{rolfsen2003knots}. 
    
    \item Stable orbits of planets require $3{+}1$ dimensions, i.e.\ three spatial dimensions and one temporal dimension, see e.g.\ \cite{ehrenfest1917way,tegmark1997dimensionality}.
    
    \item Exotic differentiable structures\footnote{We won't go to the very technical details here, but we can equip an Euclidean space with a so called differentiable structure that specifies how calculus works on that space. There is only one choice in all the Euclidean spaces except for $\R^4$, which has more.} occur only in $\R^4$ \cite{freedman1982topology,donaldson1983application}.
\end{enumerate}

\section{Mathematical bounds to almost orthogonality}
\label{sec:math_bounds_to_ao}

\begin{figure}[ht]\centering
    \includegraphics[width=1.0\textwidth]{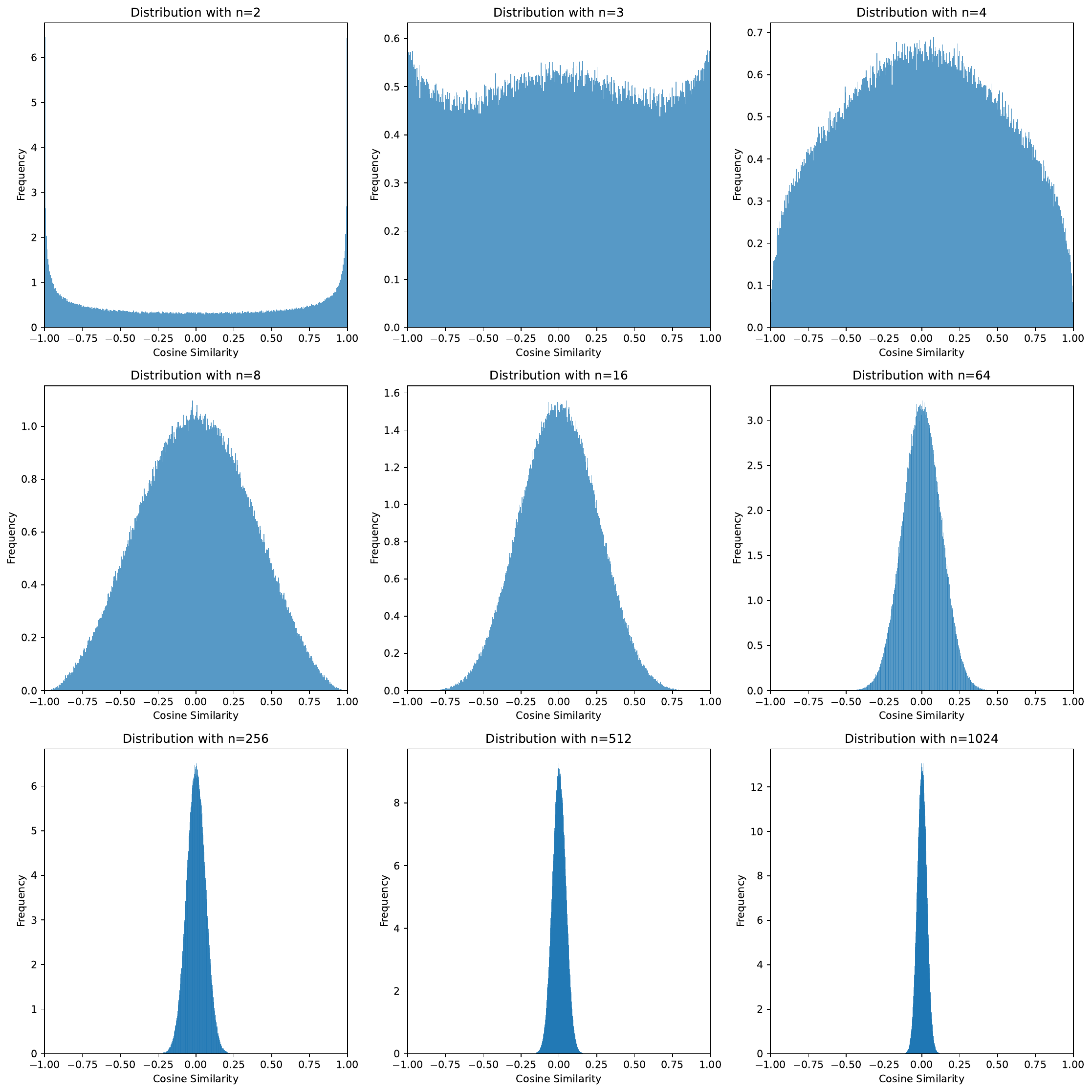}
    \caption{Cosine similarities of randomly sampled vectors in various dimensions.}
    \label{fig:RandomCosineSimilarityPlotting}
\end{figure}

Our main question in this section is to find out how many almost orthogonal vectors can we fit in a given vector space. After the previous section on how different the behavior of low and high dimensional spaces can be, we should expect the dimension to play a large role. Indeed, let's begin by doing a quick simulation. In the Figure \ref{fig:RandomCosineSimilarityPlotting} we have sampled 1000 random vectors in a few different dimensions and plotted the distribution of their cosine similarities - note that we have included for any pair of vectors $\vecv, \vecw$ both the cosine similarity $\operatorname{cs}(\vecv, \vecw)$ and $\operatorname{cs}(\vecw, \vecv)$, which does not effect the shape here as $\operatorname{cs}(\vecv, \vecw) = \operatorname{cs}(\vecw, \vecv)$.

Indeed, Vershynin \cite{vershynin2018high}\footnote{See also https://math.stackexchange.com/questions/4555166/the-probability-for-inner-product-for-two-unit-sphere-uniformly-distributed-rand.}
shows that the expected value of the cosine similarity of two random $n$-dimensional unit vectors is $1/\sqrt{n}$. This is in some way expected; see again Figure \ref{fig:sphere_dense_equators} describing how the mass of a sphere is concentrated on the equator. From this point of view if you choose one vector as your basis, then any random vector is most likely chosen from near the equator relative to the first vector, meaning that they are nearly orthogonal.

\subsection{Volume-based limits}
\label{sec:volume_based_limits}

What's the optimal amount then? How many $\eps$-almost orthogonal vectors can we have in $\rn$ for a given $\eps$ and $n$? Here it turns out that an exact answer is really hard to figure out. The problem is related to the problem of \emph{sphere packings} in high dimensions. Indeed, suppose we have a set of vectors in $\rn$ that have some upper bound on their pairwise cosine similarities. Since the cosine similarity is not altered by scalar scaling, we can assume these are all unit vectors. The upper bound on the cosine similarities implies that these vectors, when thought of as points in the $n-1$-dimensional unit sphere $\mathbb{S}^{n-1}$, have a lower bound $b = b(\eps)$ on their pairwise distances in the arc-distance of $\mathbb{S}^{n-1}$. This correspondence between the angle bound and arc-distance is one-to-one, meaning that looking for an optimal $\eps$-almost orthogonal set of vectors in $\rn$ is quantitatively \emph{equivalent} to finding a set of points $\mathbb{S}^{n-1}$ with a uniform lower bound $b$ on their pairwise distances. This in turn translates to trying to find out how many disjoint balls of radius $b$ can we have in the space $\mathbb{S}^{n-1}$. This problem is a variation of the so called \emph{sphere packing problem}, and we do not know how to solve it even in most Euclidean spaces. (The solution is currently known in dimensions 3, 8 and 24: \cite{hales2005proof,viazovska2017sphere,cohn2017sphere}.) So suffice it to say that we can't prove in this paper what is the maximum amount of $t$-almost orthogonal vectors in dimension 768.

So if an optimal bound is impossible, can we at least estimate an upper bound of these sphere packings? Let's try an approach based on volume estimates; two distinct unit vectors give rise to disjoint spherical caps, and we can calculate how many disjoint spherical caps can at most fit on a sphere based on volume. First recall Figure \ref{fig:Spherical_cap_diagram} on the terminology of spherical caps and Section \ref{sec:volumes_and_areas_of_spheres_and_balls} for the area of such caps. We want to calculate the For $\eps$-almost orthogonal vectors, the angle between them is $\alpha = \cos\inv(\eps)$, and the disjoint spherical caps have angle $\alpha/2$. 
\begin{align*}
h 
= 1 - \cos(\alpha/2) 
= 1 - \cos\left(\frac{1}{2}\cos\inv(\eps)\right).
\end{align*}

In particular, for us this boils down to
\begin{align*}
    A_h 
    =\frac{1}{2}\frac{2\pi^{n/2}}{\Gamma\left(\frac{n}{2}\right)} I_{2h-h^2} \left(\frac{n-1}{2}, \frac{1}{2} \right)
\end{align*}
and with $h = 1-\cos(\alpha/2)$ we get, by momentarily denoting $\cos(\alpha/2) \equalscolon t$, that
\begin{align*}
    2h-h^2 
    &= 2 - 2\cos(\alpha/2) - (1-\cos(\alpha/2))^2 \\ 
    &= 2 - 2t - (1-t)^2 \\
    &= 2 - 2t - 1 + 2t - t^2 \\
    &= 1-t^2 = 1 - \cos(\alpha/2)^2.
\end{align*}

So in particular we get that
\begin{align*}
    \operatorname{vol}(SC_{\alpha/2})
    = \frac{1}{2}\frac{2\pi^{n/2}}{\Gamma\left(\frac{n}{2}\right)} I_{1-\cos(\alpha/2)^2} \left(\frac{n-1}{2}, \frac{1}{2} \right).
\end{align*}

On the other hand for the area of the whole $n$-dimensional unit sphere, $\mathbb{S}^n$, we get 
\begin{align*}
    A
    = \frac{2\pi^{n/2}}{\Gamma(n/2)}.
\end{align*}
In particular, a single spherical cap takes up a fraction of
\begin{align*}
    \frac{SC_{\alpha/2}}{A} = \frac{1}{2} I_{1-\cos(\alpha/2)^2} \left(\frac{n-1}{2}, \frac{1}{2} \right)
\end{align*}
of the area of the total unit $n$-sphere. Note that each direction produces two such spherical caps at the antipodal points. Thus the inverse of twice this fraction gives an absolute area-based upper bound for the amount of $\eps$-almost orthogonal vectors in $\R^n$.

Let's get estimates for values of $\eps$ being $0$, $0.1$ and $0.01$ in the dimensions of $n=2$, $n=3$, $n=32$, $n=768$ and $n=4096$. For each $\eps$ we have $\alpha = \cos\inv(\eps)$, and for these values we wish to estimate
$$
I_{1-\cos(\alpha/2)^2} \left(\frac{n-1}{2}, \frac{1}{2} \right)^{-1},
$$
where $I_x(a,b)$ is the incomplete Beta function. The results are listed in table

\begin{table}[h]
\centering
\begin{tabular}{c|ccc}
\toprule
$n \backslash \varepsilon$ & $0.1$ & $0.01$ & $0$ \\
\midrule
$2$ & $ 2.136   $ & $ 2.013   $ & $ 2.0   $ \\
$3$ & $ 3.87   $ & $ 3.456   $ & $ 3.414   $ \\
$4$ & $ 6.605   $ & $ 5.602   $ & $ 5.504   $ \\
$8$ & $ 45.08  $ & $ 31.36  $ & $ 30.17  $ \\
$16$ & $ 1526 $ & $ 720.9  $ & $ 665.9  $ \\
$32$ & $ 1.267 \times 10^{6} $ & $ 2.784 \times 10^{5} $ & $ 2.372 \times 10^{5} $ \\
$768$ & $ 2.537 \times 10^{134} $ & $ 3.249 \times 10^{118} $ & $ 6.849 \times 10^{116} $ \\
$4096$ & $ 6.635 \times 10^{711} $ & $ 1.127 \times 10^{627} $ & $ 1.296 \times 10^{618} $ \\
\bottomrule
\end{tabular}
\caption{Some area-based bounds for the amount of $\eps$-almost orthogonal vectors in $\R^n$.}
\label{tab:inverse_beta_eps}
\end{table}

So what we note here is that even though these area estimates do give an upper bound to the almost orthogonal vectors, the estimates are quite useless. Indeed, when applied to the cosine similarity of $0$ which corresponds to actually orthogonal vectors, we get an estimate much larger than the dimension which we know to be the real bound in dimensions above three. So the \emph{geometry} of the situation is really crucial here; the dominant limit to the amount of spherical caps is not lack of volume. Indeed, as we'll discuss in Section \ref{sec:sphericalcodesandspherepacking}, the so called \emph{packing density} goes down at an exponential speed for high dimensional sphere packings.

So this purely area-based bound does tell us that we can't pack more than about $10^{100}$ of $0.1$-almost orthogonal vectors in $\R^{768}$. Later on we'll see that the actually achievable amounts seem to be in the realm of $10^{3} -- 10^{5}$, so this bound is not very useful in practice.

\subsection{The Johnson-Lindenstrauss Lemma}
\label{sec:JohnsonLindenstrauss}

We next study some existing strong results on how point clouds can be embedded into a lower ambient dimension while limiting the amount of distortion.
First we study a classical result of Johnson-Lindenstrauss that gives excellent asymptotic bounds on a slightly more complex question.

One approach to get almost orthogonal vectors in $\rn$ is to try to take a set of orthogonal vectors in a higher dimensional space, say $\R^N$, and then see if we can have a mapping $f \colon \R^N \to \rn$ that doesn't distort the cosine similarities "too much". This type of dimensionality reduction is a very important tool in e.g.\ data-analysis, compressed sensing, and other fields where we get algorithmic complexity penalties from high dimensions. A classical result to aid in this is the Johnson-Lindenstrauss lemma, see  \cite{johnson1984extensions}

\begin{lemma}[Johnson-Lindenstrauss Lemma]\label{lemma:Johnson-Lindenstrauss}
    Let $0 < \varepsilon < 1$ and let $S \colonequals \{ x_1, \ldots, x_k \} \subset \R^N$. Then 
    for any $n > \frac{8\log(k)}{\varepsilon^2}$ there exists a linear map $L\colon \R^N \to \R^n$ such that
    \begin{align}
        (1 - \varepsilon)\|x_i - x_j\|^2 
        \leq \|Lx_i - Lx_j\|^2 
        \leq (1 + \varepsilon)\|x_i - x_j\|^2
    \end{align}
    for all $x_i, x_j \in S$.
\end{lemma}

\begin{remark}
    We have used in the statement of the theorem the 
    constant 8 in the bound $n > \frac{8\log(k)}{\varepsilon^2}$, as this seems
    to be the common one listed e.g.\ in Wikipedia. However, the source for this constant are the lecture notes \cite[Lemma 2.6]{fernandezgranda2016random} which are not peer-reviewed, though they do contain a proof. The discussion on the Wikipedia articles talk page\footnote{\url{https://en.wikipedia.org/wiki/Talk:Johnson\%E2\%80\%93Lindenstrauss_lemma},viewed 24.09.2025.} lists various variants on this bound, which we summarize here:
    
    \begin{itemize}
        \item For $m>4$ we get $n > 20(\ln m)/\eps^2$ in \cite[Lemma 15.4]{mohri2018foundations}.
        \item In \cite[p.\ 300]{matousek2013lectures} they have $n>200(\ln m)/\eps^2$.
        \item \cite{dasgupta2003elementary} lists a different "style" of bound: $n > 4 \left(\varepsilon^2/2 - \varepsilon^3/3\right)^{-1}\log(k)$.
        \item The lecture notes \cite[Section 3.1.3]{duchi2024notes} prove $n > 16 (\ln m)/\eps^2$.
    \end{itemize}

    For our purposes it turns out that the exact choice here does not matter, as we shall soon see. The issue is that these results are by nature asymptotic and do not give strong estimates in the small-ish dimensions like 768. 
\end{remark}

Note that the Johnson-Lindenstrauss Lemma can be widely useful in dimension reduction. Suppose we are studying a number of grayscale images that have been taken with a megapixel resolution, i.e.\ we have a set of vectors in $\R^{1000000}$.
Many algorithms will balk at such a high dimension, e.g. if their complexity is polynomial w.r.t.\ the dimension. But if we are willing to tolerate some error, then Johnson-Lindenstrauss can help us by a lot. Indeed, imagine we have 10k images we want to study and we are okay with an error rate of $\eps = 0.1$. Now the Johnson-Lindenstrauss lemma gives us a mapping into an Euclidean space of dimension around 7000. Now suppose we have 100k images, 1M or 10M images. What we see now is that we only need to map our data to dimension of about 9k, 11k or 13k, respectively. This reduction in dimension from a million to around ten thousand can have a drastic effect in practice. Especially since the amount of images is only increasing the end dimension in a logarithmic way, as we see in these examples.

Now let's turn to see how we might use Johnson-Lindenstrauss with our problem. The natural idea would be, of course, to take a large-dimensional $\R^N$, and choose its basis plus the origin as our set of vectors. Now the first question is then about how the cosine similarity will react to the Johnson-Lindenstrauss map. First we see that
\begin{align*}
    \frac{\ip{L\vecv}{L\vecw}}{\|L\vecv\|\|L\vecw\|} 
    &= \frac{1}{2}\frac{\|L\vecv\|^2 + \|L\vecw\|^2 - \|L\vecv - L\vecw\|^2}{\|L\vecv\|\|L\vecw\|}\\
    &\leq \frac{1}{2} \frac{(1+\eps)\| \vecv\|^2 + (1+\eps)\| \vecw \|^2 - (1-\eps)\| \vecv - \vecw \|}{\sqrt{1-\eps}\|\vecv\| \sqrt{1-\eps} \| \vecw \|}.
\end{align*}
Now if $\vecv$ and $\vecw$ are orthonormal, then $\| \vecv \| = \| \vecw \| = 1$, and $\| \vecv - \vecw \|^2 = 2$, so we get that this equals to
\begin{align*}
    \frac{1}{2} & \frac{(1+\eps)\| \vecv\|^2 + (1+\eps)\| \vecw \|^2 - (1-\eps)\| \vecv - \vecw \|}{\sqrt{1-\eps}\|\vecv\| \sqrt{1-\eps} \| \vecw \|}\\
    &= \frac{1}{2} \frac{(1+\eps) + (1+\eps) - 2(1-\eps)}{(1-\eps)} \\
    &= \frac{1}{2} \frac{4\eps}{(1-\eps)} \\
    &= \frac{2\eps}{(1-\eps)}.\\
\end{align*}
Similarly we can derive a lower bound, and so we have
\begin{align}\label{eq:t_eps_relationship}
    -\frac{2\eps}{(1+\eps)}
    \leq \frac{\ip{L\vecv}{L\vecw}}{\|L\vecv\|\|L\vecw\|}
    \leq \frac{2\eps}{(1-\eps)}
\end{align}

Thus if we apply the Johnson-Lindenstrauss Lemma to the set of unit vectors in $\R^N$ appended with the origin with our
aim to generate almost orthogonal vectors in $\R^{768}$, we see that we have
to choose $\eps$ to be roughly half that of the desired threshold. For an example,
if we want the threshold to be around $\frac{1}{3}$, then we have to choose $\eps = \frac{1}{6}$.
With this choice the Johnson-Lindenstrauss lemma then gives us a bound from
\begin{align*}
    768 > 8 \frac{\ln(N+1)}{\eps^2}
    &\Leftrightarrow 96 \cdot \eps^2 > ln(N+1)\\
    &\Leftrightarrow \frac{96}{36} > \ln(N+1)\\
    &\Leftrightarrow \frac{8}{3} > \ln(N+1)\\
    &\Leftrightarrow N < e^{8/3}-1 \approx 13.4.
\end{align*}
So for the threshold of $\frac{1}{3}$, a direct application of the Johnson-Lindenstrauss technique guarantees at least 13 almost orthogonal vectors in $\R^{768}$. This is a far cry from the 768 orthogonal base vectors we know that exist in $\R^{768}$.

More generally Johnson-Lindenstrauss gives the bound
\begin{align}\label{eq:JL-dimension_bound}
    N < \exp\left(\frac{n}{8} \left(\frac{2t}{1-t}\right)^2\right)-1.
\end{align}
And so for the threshold of $0.1$-almost orthogonal vectors, to get more than the basis-guaranteed $N \geq n$ almost orthogonal vectors, then we have to have $n$ at least 29748. Though after that bound the exponential behavior takes over, and e.g.\ for $n = 40000$ and $\eps = 0.1$ we already get $N > 10^6$.

So turns out that for our "not quite so large" dimension of 768, the Johnson-Lindenstrauss result does not yield results applicable for our almost orthogonal vectors. The issue here stems largely from the fact that the Johnson-Lindenstrauss result is an asymptotic result and gives a much stronger output than just almost orthogonal vectors - it roughly preserves distances of an \emph{arbitrary} set of vectors.

We know, however, that more almost orthogonal vectors can be fitted into $\R^{768}$. For example, the vectors $(1,1,1,\ldots)/\sqrt{768}$ and $(1,-1,1,-1,\ldots)/\sqrt{768}$ have a very small cosine similarity (below $0.04$) both with all the standard basis vectors and each other. With a bit of clever combinatorics it's not hard to get many more such examples. Since these theoretical results don't seem to provide impressive results in this realm of dimensions, let us turn to less theoretical approaches. Though as we will later see the underlying idea behind the proof of the Johnson-Lindenstrauss Lemma is still very applicable for our purposes.

\subsection{Packing density, spherical codes and the Ka\-ba\-ti\-ans\-kii-Le\-ven\-sh\-tein bound}
\label{sec:sphericalcodesandspherepacking}

In this section we'll briefly mention a few related topics. We won't dive deeper to these ideas, but note them for context.

As we noted in Section \ref{sec:volume_based_limits}, the pure volume-based bounds on almost orthogonality were completely useless in high dimensions. The issue here is related to the so called \emph{packing density} of spherical packings. We refer to \cite{cohn2016conceptual} for an overview of the topic of sphere packing and packing density, but the main point for us is that for high dimensions sphere packings are not dense; they cover an exponentially diminishing fraction of the ambient space as the dimension increases. We saw a part of this effect on the rightmost column of Table \ref{tab:inverse_beta_eps} -- there the true maximum amount of disjoint spherical caps with antipodal pairs is dictated by linear algebra, but the bounds provided by essentially area density grow in some superexponential fashion.

So almost orthogonal vectors correspond to pairwise disjoint spherical caps. Another way of phrasing this is that they correspond to a set of points on a sphere with an upper bound to their pairwise inner products. Such point sets on a sphere are known as \emph{spherical codes} and have direct applications in coding theory. We refer the reader to \cite{cohn2014sphere} for details, but note that the study of spherical codes is also an unsolved field of study.

A related problem is known as \emph{the Tammes problem} (see e.g.\ \cite{musin2018five}) where we ask to minimize the pairwise distances of a set of points on the sphere. Anecdotally we've seen the Tammes problem used as a term in the realm of "small balls and distances", whereas the "almost orthogonal" approach is more in the realm of large distances as the points induced on the sphere aim to be almost a quarter apart. 

Note again that both the Tammes problem and spherical codes usually focus on \emph{points on a sphere} while we work on \emph{directions in space}. Any direction naturally induces two antipodal points on a sphere, but this extra assumption of antipodal inclusion is not present on the more general questions of spherical codes or the Tammes problem.

\section{Simulation approaches}

Our aim is to next turn into generating collections of almost orthogonal vectors in various Euclidean spaces. As we've learned so far, there is very little hope of proving anything about the optimality of any of these collections. But if different methods seem to yield similar results, we could consider that this as weak evidence that an actual mathematical limit might be near.

\subsection{Methods}

Our three main generation approaches are:
\begin{enumerate}
    \item Sampling random unit vectors. 
    \item Random projections of higher-dimensional bases in the spirit of the proof idea behind the Johnson-Lindenstrauss Lemma.
    \item Energy minimization through ML methods.
\end{enumerate}

Each of these have their own subsection below, explaining the idea of the approach. Besides these generative techniques, we also test \emph{pruning} methods, i.e.\ generating extra vectors and then dropping a subset with an aim of reducing the maximum absolute cosine similarity.

\subsubsection{Randomly sampled unit vectors}

For us this is the most straightforward method. If we sample unit vectors randomly, then they should not have massive pairwise cosine similarities. Indeed, the expected value of the cosine similarity between two randomly sampled unit vectors in $\R^n$ should be around $\frac{1}{\sqrt{n}}$. (See again the discussion in Section \ref{sec:math_bounds_to_ao}.)

It is not trivial\footnote{Or, it was not trivial for us.} to derive what the expected maximum absolute cosine similarity then would be for $N$ randomly generated unit vectors, which is why we turned to simulations. Though see again Figure \ref{fig:RandomCosineSimilarityPlotting} where we show a descriptive example of the distributions of pairwise cosine similarities of random unit vectors in various dimensions.

After some heuristical tests, we settled to generating twice the amount of vectors here and then pruning half of them.

\subsubsection{Random projections of higher-dimensional bases}

Instead of raw sampling, we can apply the procedure behind the proof of the Johnson-Lin\-den\-strauss Lemma discussed earlier (Section \ref{sec:JohnsonLindenstrauss}). Here the basic idea is that if we take a collection of vectors in a high dimension, then its projection to a random subspace has a non-zero probability to weakly preserve its various (geo)metric properties. So for our purposes of trying to find almost orthogonal collections of vectors, a natural approach here would be to take a random orthonormal basis in a high-dimensional space $\R^N$ and project it to a lower $n$-dimensional space subspace $V \subset \R^N$. 

Due to rotational symmetry, it suffices to either take the standard basis of $\R^N$ and project it to a random subspace or to take a random orthonormal basis and project it to the subspace spanned by the first $n$ coordinates i.e.\ take the first $n$ coordinates of the vectors. We'll opt for the latter here, though the former could probably be optimized to run a bit faster and with less memory.

Also in this generation approach we opted to generating twice the amount of vectors and pruning half.

\subsubsection{Minimizing the absolute value of cosine similarity}

Partially motivated by a video of the math Youtuber 3Blue1Brown\footnote{\url{https://www.youtube.com/watch?v=9-Jl0dxWQs8}} we wanted to try out a more machine-learning based method. Here we create a nice collection of vectors by minimizing a type of energy function on the pairwise distances of vectors. To this end we take a sequence of unit vectors representing directions $[v_1, \ldots, v_{n}]$, append them with their antipodal vectors $[-v_1, \ldots, -v_{n}]$, and then look at the energy
\begin{align*}
    \sum_{i<j} \frac{1}{(d(w_i, w_j))^p},
\end{align*}
where the $(w_j)_{j=1}^{2n}$ is the set of vectors and their antipodal vectors. We start with a random sampling of vectors, and then "nudge" each of the vectors $v_i$ in direction that most reduces this energy.

In this approach the pruning methods proved, unsurprisingly, mostly detrimental.

\subsection{Combinations and hyperparameters}

We tested several variations and combinations of these generation and pruning techniques. For the energy minimization there was little effect on the results based on whether we started with a random collection of vectors or e.g.\ something generated with the random projections method. We also tested some variations on how we sample random unit vectors, but the results were indistinguishable.

For the energy minimization there were various hyperparameters, but as our aim here is to show more qualitative than quantitative results we won't go to details on our choices. The source will be made available. We will, however, mention that the choice for th exponent in our energy functional had a big effect on how fast the process converged. Figures \ref{fig:steps_in_optimization_512} and \ref{fig:steps_in_optimization_2048} below show how the maximum absolute cosine similarity evolves with various exponent choices. We see that the evolution converges much faster when the exponent is increased until around exponent 16 after which results start to worsen. This is most likely due to the fact that the higher exponent produces stronger gradients, until the floating point numbers start to round off to zero. There is probably some clever way of optimizing this procedure, but we leave it out of the scope of this article.

\begin{figure}[ht]\centering
    \includegraphics[width=0.8\textwidth]{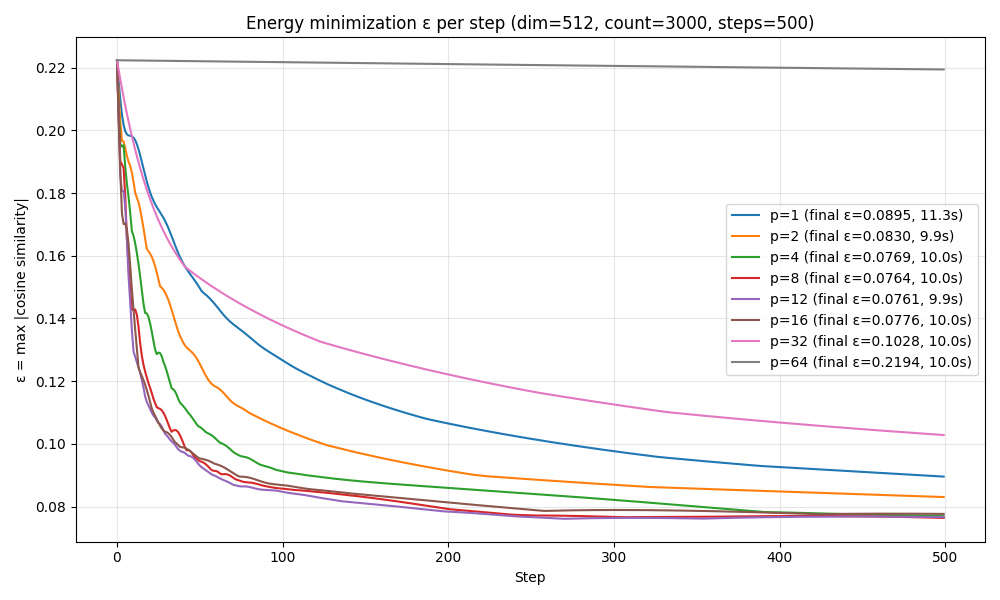}
    \caption{The effect of the energy exponent and step count in generating 3000 unit vectors in $\R^{512}$ minimizing their energy.}
    \label{fig:steps_in_optimization_512}
\end{figure}

\begin{figure}[ht]\centering
    \includegraphics[width=0.8\textwidth]{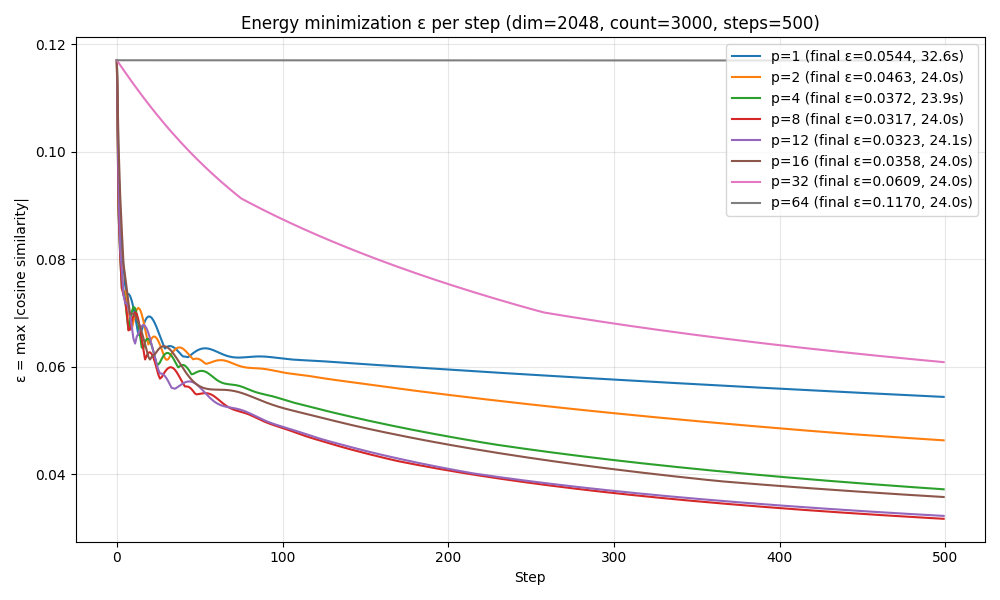}
    \caption{The effect of the energy exponent and step count in generating 3000 unit vectors in $\R^{2048}$ and minimizing their energy.}
    \label{fig:steps_in_optimization_2048}
\end{figure}

We also tested a few variants of the energy minimization approach where we used other loss functions that targeted the $\max$ or top $N$ absolute cosine similarities directly. They did not prove to be very effective in practice.\footnote{We also note that the 3Blue1Brown youtube video mentioned earlier also used a gradient flow method but aimed to minimize the \emph{average} cosine similarity which produced good average results but strong outliers.}

In the pruning methods we tried a few variations of the amount of oversampling and then pruning. The factor of $2.0$ was chosen quite heuristically as a balance between compute time and results. 
We note that besides just dropping the "worst offenders" from the generated set one by one, a more subtle approach could be to try and detect what is in some sense the best subset of vectors from a given generated set. If we phrase this as a graph problem where each vector is a node and two nodes are connected when the corresponding two vectors have their absolute cosine similarity below some threshold, then this turns into the NP-complete \emph{Clique problem}. This formulation does naturally discard a lot of geometric information, and finding a best (or good) almost orthogonal subset might be computationally much easier in our case. But with the problem adjacent to an NP-complete problem we would be wary.

\subsection{Generation results}

\begin{figure}[h!]
    \centering
    \includegraphics[width=1.0\textwidth]{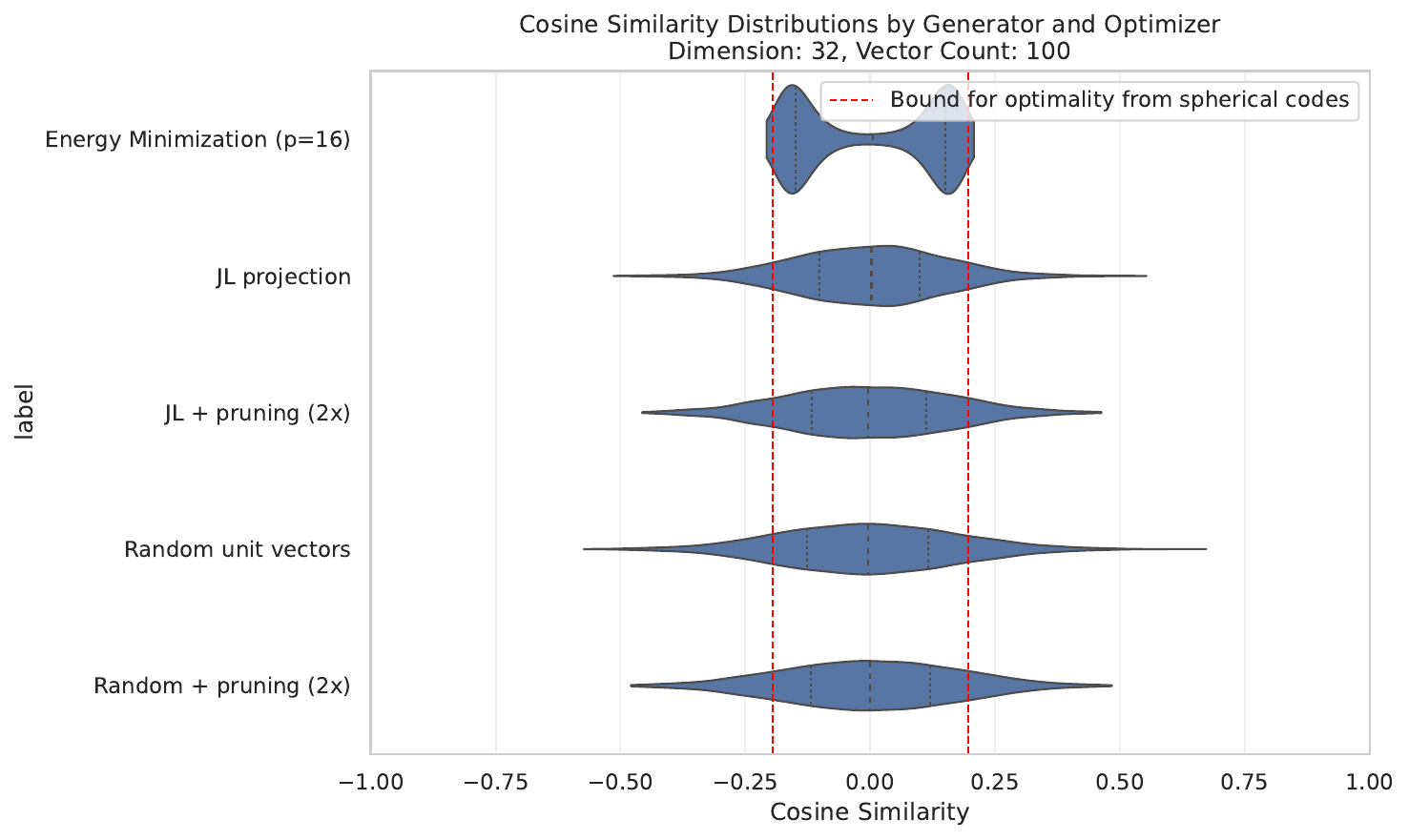}
    \caption{Violin plot of the pairwise cosine similarity distributions generated by various methods in dimension $32$ and $100$ vectors. The "known optimal" is from spherical codes and most likely unattainable with directions.}
    \label{fig:violin_comp_32_100}
\end{figure}

\begin{figure}[h!]
    \centering
    \includegraphics[width=1.0\textwidth]{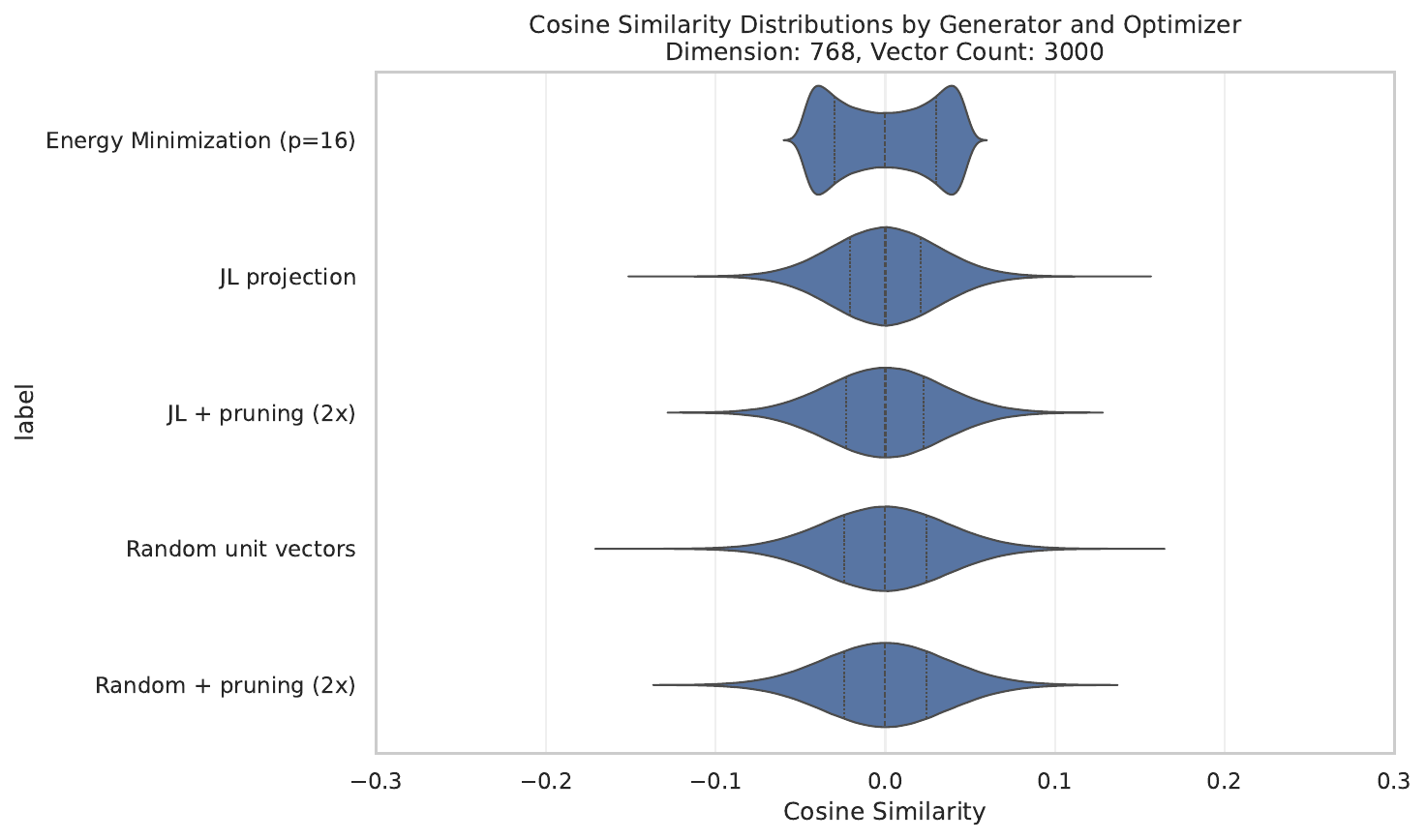}
    \caption{Violin plot of the pairwise cosine similarity distributions generated by various methods in dimension $768$ and $3,000$ vectors.}
    \label{fig:violin_comp_768_3k}
\end{figure}

\begin{figure}[h!]
    \centering
    \includegraphics[width=1.0\textwidth]{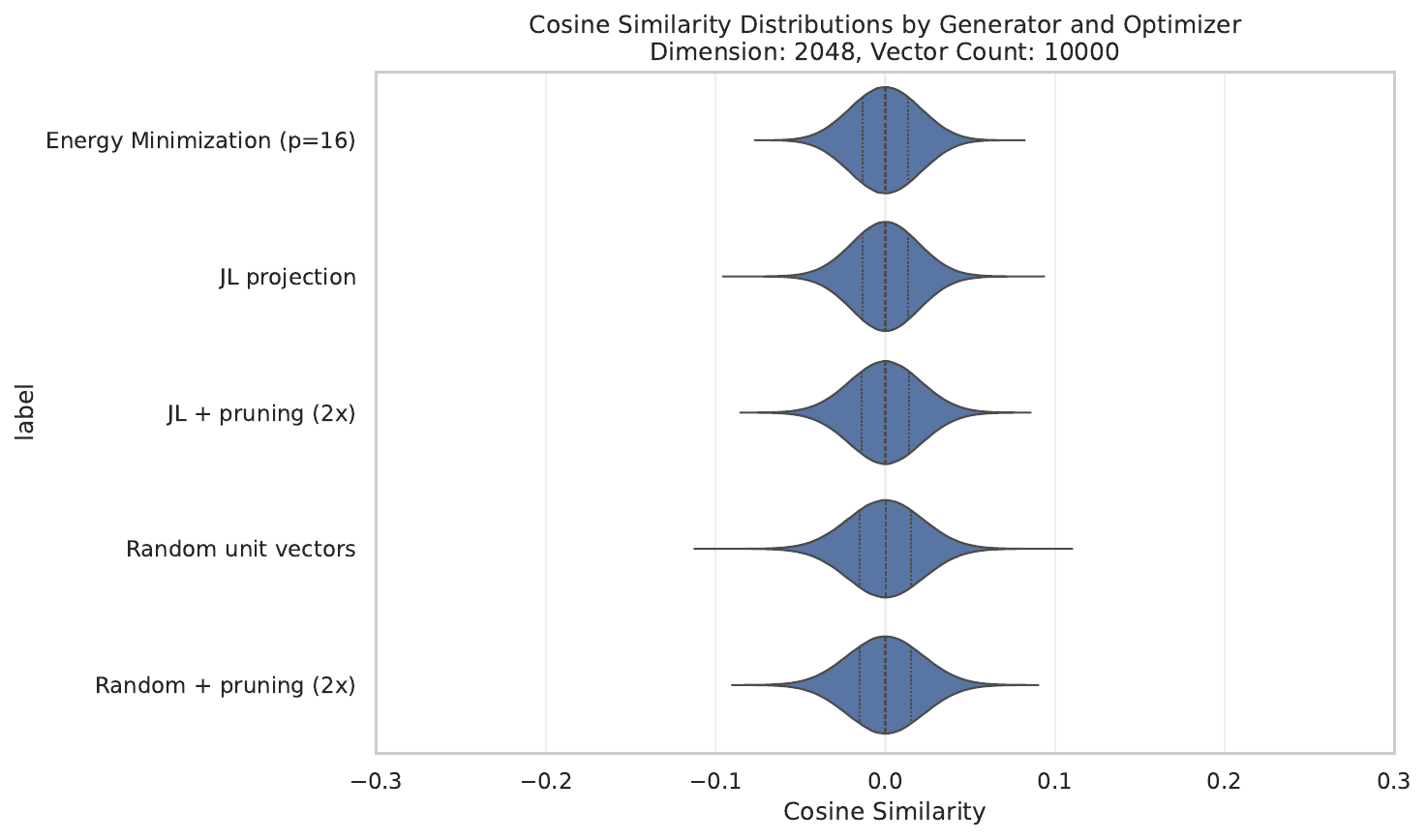}
    \caption{Violin plot of the pairwise cosine similarity distributions generated by various methods in dimension $2048$ and $10,000$ vectors.}
    \label{fig:violin_comp_2048_10k}
\end{figure}

First of all, in Figures \ref{fig:violin_comp_32_100}, \ref{fig:violin_comp_768_3k}, and \ref{fig:violin_comp_2048_10k} we show what kind of distribution of cosine similarities our generators produce in dimensions 32, 768 and 2048 for specific sizes of vector collections. What we note is that, especially in the lower dimensions, the energy minimizers create a bimodal distributions. The others seem more normal distribution-y. 

\begin{figure}[h!]
    \centering
    \includegraphics[width=0.9\linewidth]{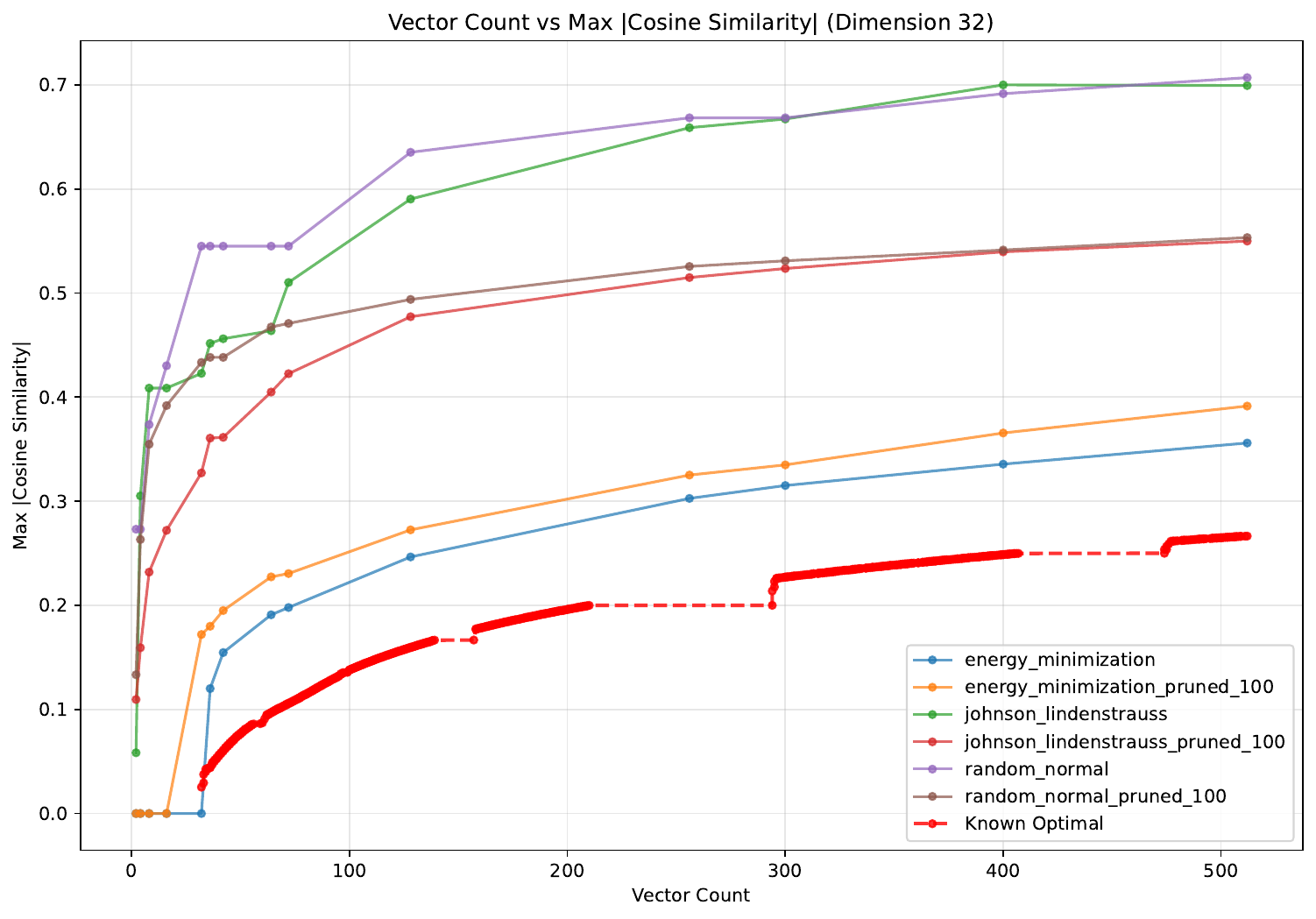}
    \caption{All vector counts in dimension 32.}
    \label{fig:comparison_graph_dim_32_all}
\end{figure}

\begin{figure}[h!]
    \centering
    \includegraphics[width=0.9\linewidth]{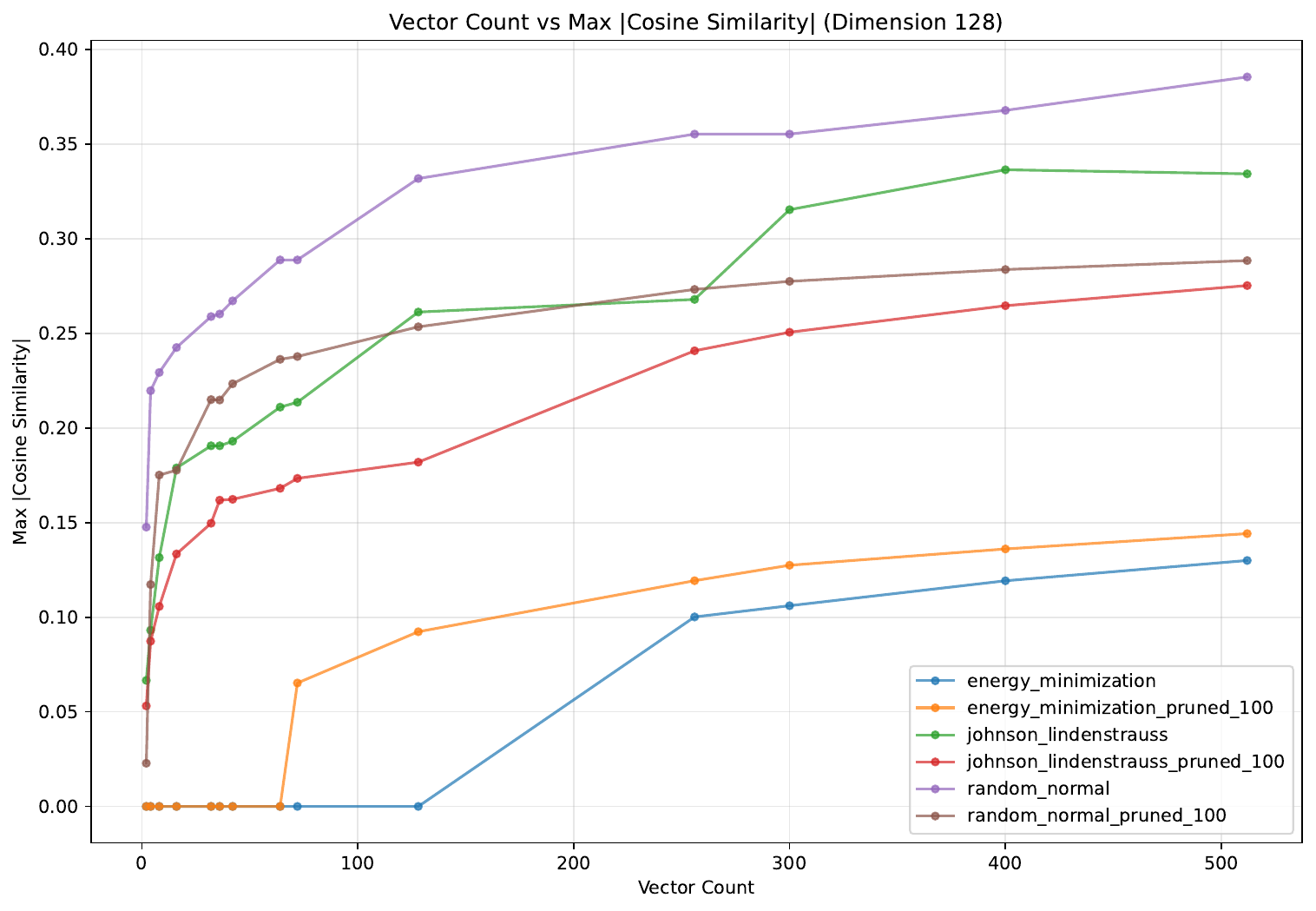}
    \caption{All vector counts in dimension 128.}
    \label{fig:comparison_graph_dim_128_all}
\end{figure}

\begin{figure}[h!]
    \centering
    \includegraphics[width=0.9\linewidth]{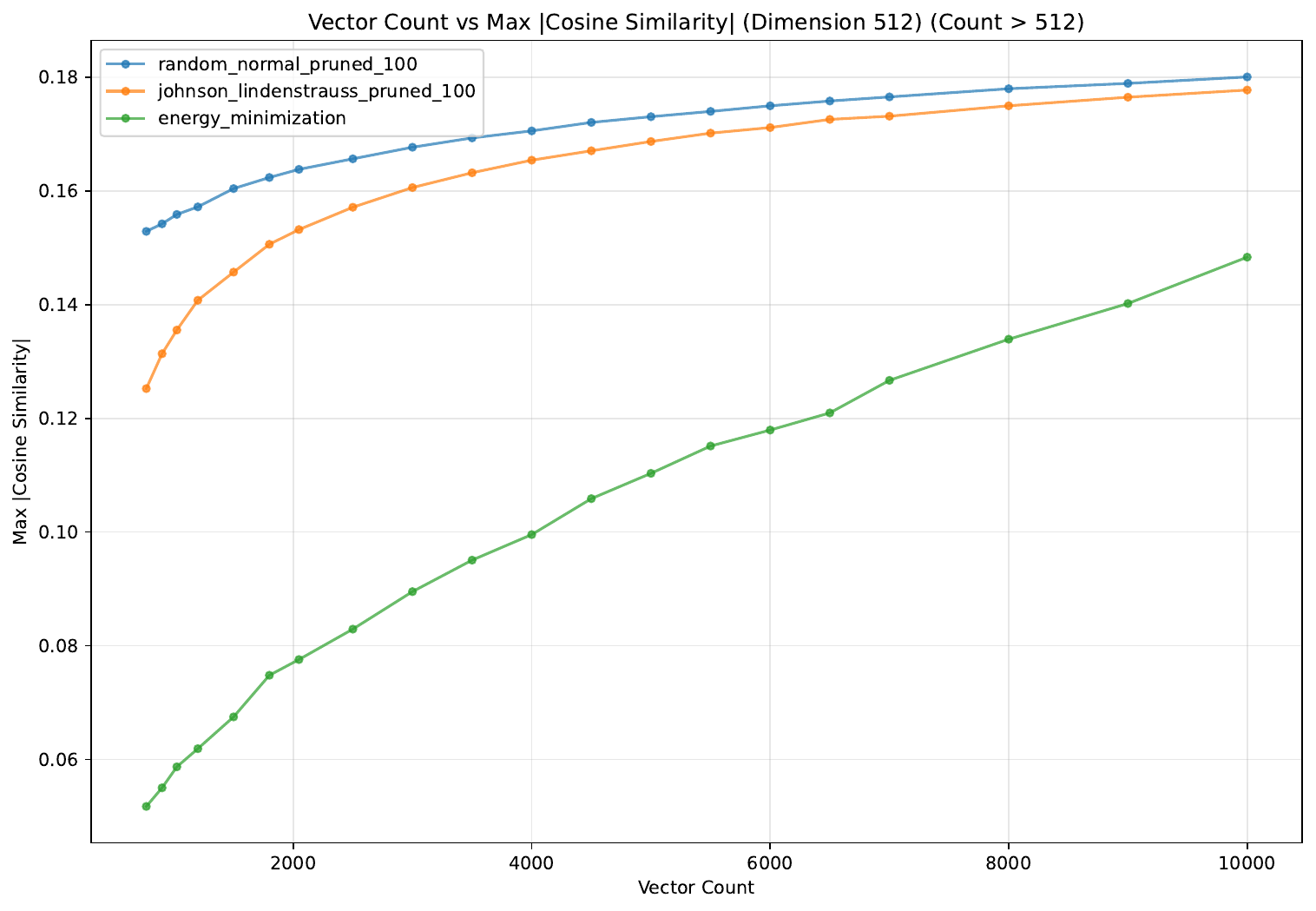}
    \caption{Generating more vectors than the dimension in dimension 512.}
    \label{fig:comparison_graph_dim_512_filtered}
\end{figure}

\begin{figure}[h!]
    \centering
    \includegraphics[width=0.9\linewidth]{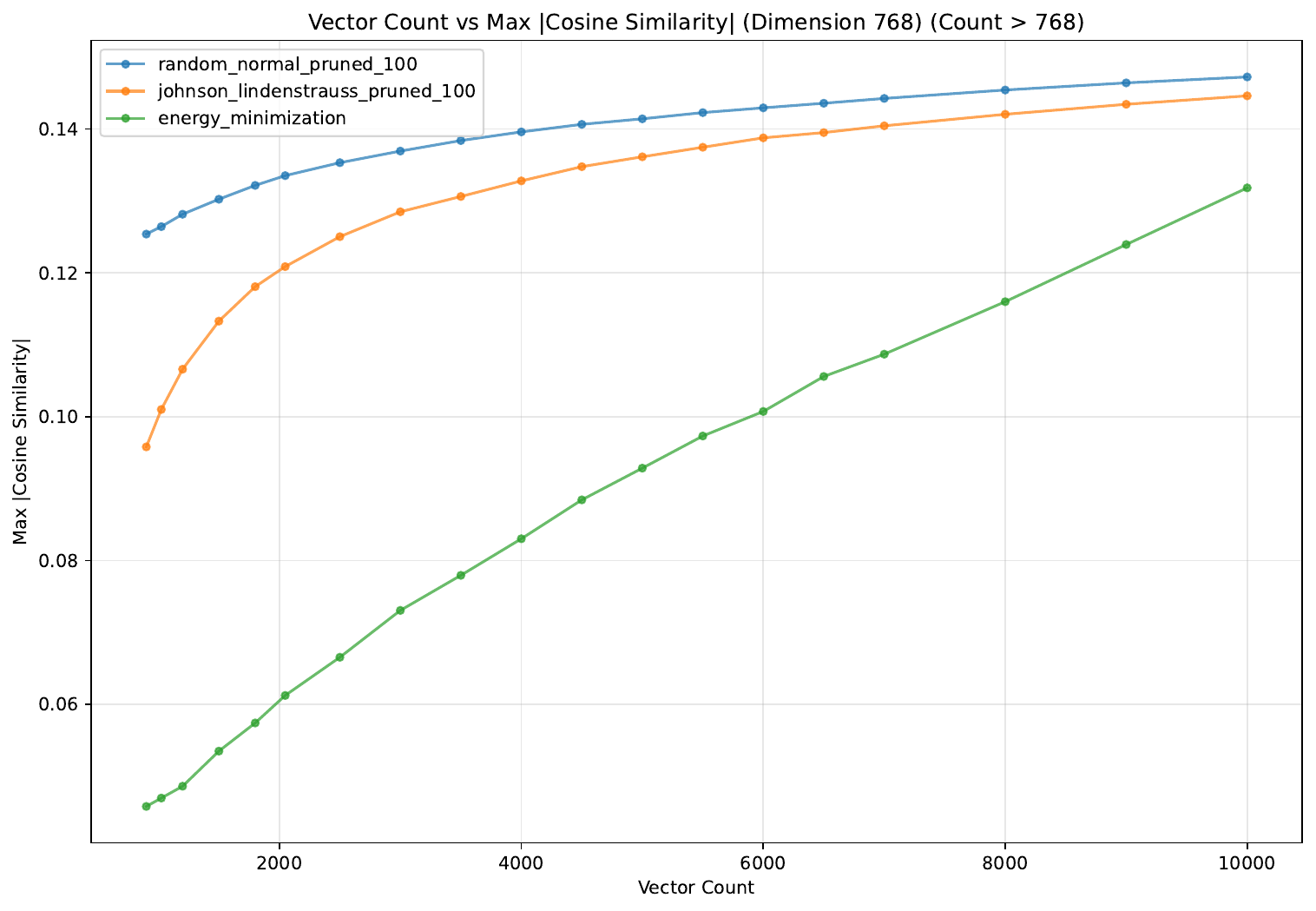}
    \caption{Generating more vectors than the dimension in dimension 768.}
    \label{fig:comparison_graph_dim_768_filtered}
\end{figure}

\begin{figure}[h!]
    \centering
    \includegraphics[width=0.9\linewidth]{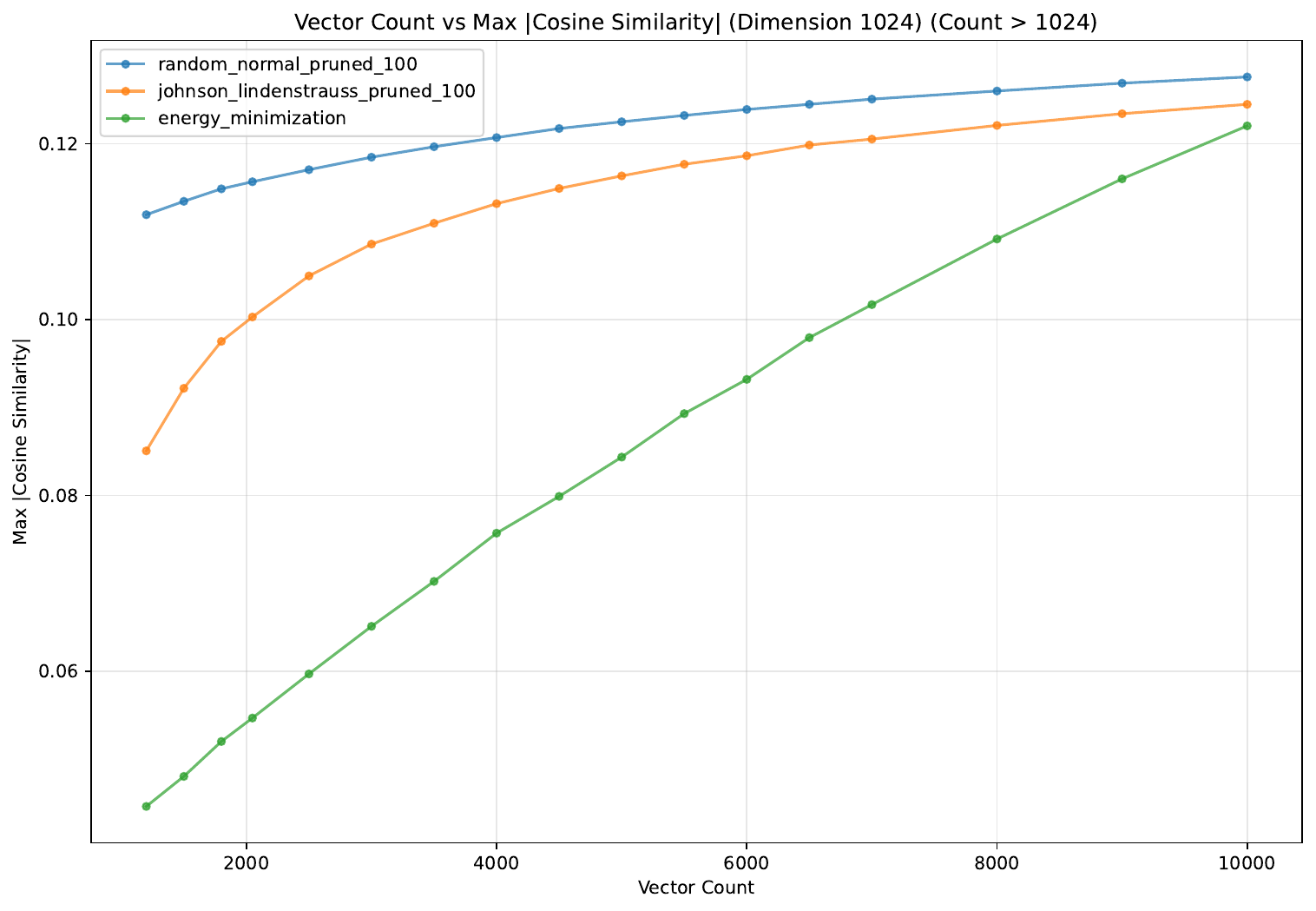}
    \caption{Generating more vectors than the dimension in dimension 1024.}
    \label{fig:comparison_graph_dim_1024_filtered}
\end{figure}

Then in Figures 
\ref{fig:comparison_graph_dim_32_all}, 
\ref{fig:comparison_graph_dim_128_all}, 
\ref{fig:comparison_graph_dim_512_filtered}, 
\ref{fig:comparison_graph_dim_768_filtered}, and
\ref{fig:comparison_graph_dim_1024_filtered}
we show how the various techniques compare with a varying amount of vectors in dimensions 32, 128, 512, 768 and 1024. In dimension 32 we also have plotted some data from \cite{cohn2023spherical} that shows the max cosine similarity of the best known corresponding sphere packing.\footnote{Here we again emphasize that our collections of almost orthogonal vectors correspond to collections of pairs of antipodal points on the sphere. The spherical packing minimizes have no such limitation of needing to have the antipodal point also present, and thus can naturally yield better results.}

\subsection{Simulation conclusions}

In the following Table \ref{table:conclusions} we've collected our best results and the method we used to acquire it. As thoroughly discussed above, these should not be considered to exact or even near exact but at least they should be indicative on what simple approaches can produce.

\begin{table}[h!]
\begin{tabular}{r|llllllll}
\hline
\textbf{$k \backslash d$} 
     & \textbf{32}  & \textbf{64}  & \textbf{128}  & \textbf{256}  
     & \textbf{512} & \textbf{768} & \textbf{1024} & \textbf{2048} \\
\hline
40   & 0.1184  & 0  & 0  & 0  & 0  & 0  & 0  & 0 \\
60   & 0.1637  & 0  & 0  & 0  & 0  & 0  & 0  & 0 \\
100   & 0.2092  & 0.1099  & 0  & 0  & 0  & 0  & 0  & 0 \\
200   & 0.2577  & 0.1503  & 0.0788  & 0  & 0  & 0  & 0  & 0 \\
400   & 0.3110  & 0.1903  & 0.1102  & 0.0564  & 0  & 0  & 0  & 0 \\
600   & 0.3416  & 0.2129  & 0.1258  & 0.0706  & 0.0330  & 0  & 0  & 0 \\
800   & 0.3713  & 0.2286  & 0.1388  & 0.0786  & 0.0422  & 0.0252  & 0  & 0 \\
1000   & 0.3887  & 0.2407  & 0.1482  & 0.0865  & 0.0472  & 0.0320  & 0  & 0 \\
1500   & 0.4109  & 0.2654  & 0.1654  & 0.1008  & 0.0573  & 0.0410  & 0.0323  & 0 \\
2000   & 0.4241  & 0.2838  & 0.1764  & 0.1075  & 0.0644  & 0.0467  & 0.0376  & 0 \\
3000   & 0.4439  & 0.3116  & 0.1968  & 0.1218  & 0.0749  & 0.0567  & 0.0458  & 0.0316 \\
5000   & 0.4752  & 0.3432  & 0.2219  & 0.1406  & 0.0900  & 0.0673  & 0.0567  & 0.0406 \\
8000   & 0.4976  & 0.3658  & 0.2455  & 0.1604  & 0.1034  & 0.0823  & 0.0707  & 0.0525 \\
     \hline
\end{tabular}
\caption{Best achieved $\max |\cos|$ by $(k, d)$.}
\label{table:conclusions}
\end{table}

\section{Afterword}\label{sec:afterword}

So after all of our analysis we can say that there is clearly an interesting phenomenon present in finding sets of almost orthogonal vectors in high dimensions. The existing theoretical bounds are still quite far from practical results.

From the point of view of e.g.\ the BERT language model, what we've noted is that with thresholds around $0.1$ we can fit an order of magnitude more almost orthogonal vectors in $\R^{768}$ than what we get from the basis vectors, which goes to show tha there is "a lot of space in the geometry" of $\R^{768}$ for lossily storing information as directions.

Many of the approaches mentioned here would benefit from a deep dive analysis, and we list some of them in the next subsection.

\subsection{Possible future avenues of study}

We list here some research questions that we found interesting, but did not have time to go more into and are not planning to work on in the near future.
We remind the reader that a part of the codebase used for the simulations in this work is available at \url{https://github.com/ramiluisto/AlmostOrthogonalGenerators}.

\begin{enumerate}
    \item Let's assume that the internal model of an language model has fixed $K$ content directions in $\R^{768}$. Then look at various shapes in $\R^K$, and map these to $\R^{768}$ in some way that maps the $K$ basis vectors to $K$-almost orthogonal vectors. (E.g.\ Johnson-Lindenstrauss style random projection techniques). Look at the distributions of the cosine similarities, is it common for the distributions to resemble those in Figure \ref{fig:violin_plot_for_model_cosims}?
    \item Compressed sensing (see e.g.\ \cite{Mackenzie2009Pixel} or \cite{donoho2006compressed}) feels to have some similar ideas of doing high-dimensional analysis in smaller dimensions. In a few brief attempts we didn't find a good connection yet, theoretically or numerically, but we feel there might be something useful in here.
    \item A deeper dive to the ideas behind the Kabatianskii-Levenshtein bound would be interesting. We doubt that the bounds are closer to optimal than our numerical examples, but it would be valuable to know.
\end{enumerate}

\newcommand{\etalchar}[1]{$^{#1}$}

\end{document}